\documentclass[11pt,reqno]{amsart}

\usepackage{amsmath,amssymb,mathrsfs}
\usepackage{graphicx,cite,times}
\usepackage{cases}
\newtheorem{theo}{Theorem}[section]

\newtheorem{lemm}[theo]{Lemma}

\numberwithin{equation}{section}

\begin{document}

\title[Time-Domain Maxwell's Equation in an Unbounded Structure]{Electromagnetic
Scattering for Time-Domain Maxwell's Equations in an Unbounded Structure}

\author{Yixian Gao}
\address{School of Mathematics and Statistics, Center for Mathematics and
Interdisciplinary Sciences, Northeast Normal University, Changchun, Jilin 130024, P.R.China}
\email{gaoyx643@nenu.edu.cn}

\author{Peijun Li}
\address{Department of Mathematics, Purdue University, West Lafayette, IN 47907,
USA.}
\email{lipeijun@math.purdue.edu}

\thanks{The research of YG was supported in part by NSFC grant 11571065 and
Jilin science and technology development project. The research of PL was
supported in part by the NSF grant DMS-1151308.}

\keywords{Time-domain Maxwell's equations, unbounded rough surfaces, Laplace
transform, stability, a priori estimates}

\begin{abstract}
The goal of this work is to study the electromagnetic scattering problem of
time-domain Maxwell's equations in an unbounded structure. An exact transparent
boundary condition is developed to reformulate the scattering problem into an
initial-boundary value problem in an infinite rectangular slab.
The well-posedness and stability are established for the reduced problem. Our
proof is based on the method of energy, the Lax--Milgram lemma, and the
inversion theorem of the Laplace transform. Moreover, a priori estimates with
explicit dependence on the time are achieved for the electric field by directly
studying the time-domain Maxwell equations.
\end{abstract}

\maketitle

\section{Introduction}

Consider the propagation of an electromagnetic wave which is excited by electric
current density and is scattered by infinite rough surfaces. An infinite rough
surface is a non-local perturbation of an infinite plane surface such that the
whole surface lies within a finite distance of the original plane. The goal of
this paper is to examine the electromagnetic scattering problem of time-domain
Maxwell's equation in such an unbounded structure. The problem studied in this
work falls into the class of rough surface scattering problems, which arise from
various applications such as modeling acoustic and electromagnetic wave
propagation over outdoor ground and sea surfaces, optical scattering from the
surface of materials in near-field optics or nano-optics, detection of
underwater mines, especially those buried in soft sediments. These problems are
widely studied in the literature and various methods have been
investigated \cite{Ogilvy1991wave, Saillard2001rigorous, Voronovich1994wave,
Warnick2001numerical, Desanto2002scattering, Elfouhaily2004critical}. 

The infinite rough surfaces scattering problems are quite challenging due to
unbounded domains. The usual Sommerfeld (for acoustic waves) or
Silver--M\"{u}ller (for electromagnetic waves) radiation condition is not valid
any more \cite{Zhang1998accoustic, Arens2005radiation}. The Fredholm alternative
theorem is not applicable due to the lack of compactness result. We refer to
\cite{Chandler2006acoustic, Chandler2005existence, Chandler1999scattering,
Li2012analysis, Lechleiter2010variational} for some
mathematical studies on the two-dimensional Helmholtz equation. The rigorous
mathematical analysis is very rare for the three-dimensional Maxwell equations.
In \cite{Li2011electromagnetic}, the electromagnetic scattering by unbounded
rough surfaces was considered by assuming that the medium was lossy in
the entire space. The well-posedness was established by a direct application of
the Lax--Milgram theorem after showing that the sesquilinear form was coercive.
In \cite{Haddar2011electromagnetice}, the authors considered the
electromagnetic scattering by an unbounded dielectric medium which was
deposited on a perfectly electrically conducting plate. Based on the limiting
absorption principle, the problem was shown to have a unique weak solution from
a prior estimates. The magnetic permeability was assumed to be a constant and
the electric current was assumed to be divergence free. The assumption was also 
restrictive for the dielectric permittivity. In \cite{LiZhengZhengMMAS}, the
generalized Lax--Milgram theorem was adopted to establish the well-posedness
for the same problem as that in \cite{Haddar2011electromagnetice}. Although
all the assumptions were relaxed, such as the magnetic permeability was allowed
to be a variable function and the divergence free condition was removed for the
electric current, the assumption was still quite restrictive for the dielectric
permittivity. Despite the tremendous effort made so far, it is still unclear
what the least restrictive conditions are for the dielectric permittivity and
the magnetic permeability to assure the well-posedness of the time-harmonic
Maxwell equations in unbounded structures. Ultimately, one wishes to answer the
following question: Is the scattering problem in unbounded structures well-posed
for the real and dielectric permittivity and magnetic permeability?

In this work, an initial attempt is made to study the time-domain
electromagnetic scattering by infinite rough surfaces for the most difficult
case of the time-harmonic counterpart: the dielectric permittivity and the
magnetic permeability are assumed to be real and bounded measurable functions.
An exact time-domain transparent boundary condition (TBC) is developed to reduce
the scattering problem into an initial-boundary value problem in an infinite
rectangular slab. To show the well-posedness, we split the reduced problem into
two sub-problems: one has homogeneous initial conditions and another has a
homogeneous boundary condition. Hence two auxiliary scattering problems need to
be considered: one is the time-harmonic Maxwell equations with a complex
wavenumber and another is the time-domain Maxwell equations with perfectly
electrically conducting (PEC) boundary condition. Based on the stability
results for the auxiliary problems, the reduced problem is shown to have a
unique solution. Our proofs rely on the Laplace transform, the Lax--Milgram
theorem, and the Parseval identity between the frequency domain and the
time-domain. Moreover, a priori estimates, featuring an explicit dependence on
the time and a minimum regularity requirement of the initial conditions and the
source term, are established for the electric field by studying directly the
time-domain Maxwell equations. 

The time-domain scattering problems have recently attracted considerable
attention due to their capability of capturing wide-band signals and modeling
more general material and nonlinearity \cite{ChenMonk2014, Jin2009, JiHuang2013,
Riley2008, WangWang2012}, which motivates us to tune our focus from seeking the
best possible conditions for those physical parameters to the time-domain
problem. Comparing with the time-harmonic problems, the time-domain problems
are less studied due to the additional challenge of the temporal dependence.
The analysis can be found in \cite{WangWang2014, Chen2008Maxwell} for the
time-domain acoustic and electromagnetic obstacle scattering problems. We refer
to \cite{LiWangWood2015} for the analysis of the time-dependent electromagnetic
scattering from a three-dimensional open cavity. Numerical solutions can be
found in \cite{Li2015, Veysoglu1993} for the time-dependent wave scattering by
periodic structures.

The paper is organized as follows. In section \ref{PF}, the model problem is 
introduced and reduced equivalently into an initial-boundary value problem by
using a TBC. Some regularity properties of the trace operator are presented.
In section \ref{at}, two auxiliary problems of Maxwell's equations are
discussed to pave the way for the analysis of the main result in section \ref
{rp}. Section \ref{rp} is devoted to the well-posedness and stability of the
reduced time-domain Maxwell equations and a priori estimates of the solution.
The paper is concluded with some general remarks in section \ref{rem}.

\section{Problem formulation}\label{PF}

In this section, we introduce the model problem and present an exact time-domain
transparent boundary condition to reduce the scattering problem into
an initial-boundary value problem in an infinite rectangular slab. 

\subsection{A model problem}

\begin{figure}
\centering
\includegraphics[width=0.3\textwidth]{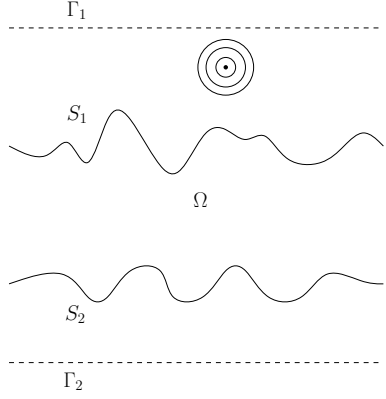}
\caption{Problem geometry of the electromagnetic scattering by an unbounded
structure}
\label{pg}
\end{figure}

Let us first introduce the problem geometry which is shown in Figure \ref{pg}.
Let $S_j, j=1, 2$ be Lipschitz continuous surfaces which are embedded in the
infinite rectangular slab
\[
\Omega=\{\boldsymbol x=(x, y, z)^\top\in \mathbb R^3: h_2 < z < h_1\},
\]
where $h_j$ are constants. Denote by $\Gamma_j=\{\boldsymbol{x}: z=h_j\}$ the
two plane surfaces which enclose $\Omega$. Let $\Omega_1=\{\boldsymbol
x: z>h_1\}$ and $\Omega_2 =\{\boldsymbol x: z <h_2\}$. The medium is assumed to
be homogeneous in $\Omega_j$, but it is allowed to be inhomogeneous in $\Omega$.

The electromagnetic field is governed by the time-domain Maxwell equations in
$\mathbb R^3$ for $t>0$:
\begin{equation}\label{EH}
 \nabla \times \boldsymbol E (\boldsymbol x, t)+\mu \partial_t \boldsymbol H
(\boldsymbol x, t)=0,\quad \nabla \times \boldsymbol H (\boldsymbol x, t)-
\varepsilon \partial_t \boldsymbol E (\boldsymbol x, t)=\boldsymbol J
(\boldsymbol x, t),
\end{equation}
where $\boldsymbol E$ is the electric field, $\boldsymbol H$ is the magnetic
field, $\boldsymbol J$ is the electric current density which is assumed to be
compactly supported in $\Omega$, the material parameters $\varepsilon$ and $\mu$
are the dielectric permittivity and the magnetic permeability, respectively. We
assume that $\varepsilon \in L^{\infty} (\mathbb R^3)$ and $ \mu \in L^{\infty}
(\mathbb R^3)$ satisfy
\[
0 < \varepsilon_{\min} \leq \varepsilon \leq \varepsilon_{ \max} <\infty, \quad
0 < \mu_{\min} \leq  \mu \leq \mu_{\max} < \infty,
\]
where $\varepsilon_{\min}, \varepsilon_{\max}, \mu_{\min}, \mu_{\max}$ are
constants. Since the medium is homogeneous in $\Omega_j$, there exist constants
$\varepsilon_j$ and $\mu_j$ such that
\[
\varepsilon (\boldsymbol x)= \varepsilon_j, \quad \mu (\boldsymbol x)=\mu_j
\quad {\rm in} ~\Omega_j.
\]
The system is constrained by the initial conditions:
\begin{equation}\label{IC}
\boldsymbol E|_{t=0}=\boldsymbol E_0, \quad
\boldsymbol H|_{t=0}=\boldsymbol H_0\quad  {\rm in}~\mathbb R^3,
\end{equation}
where $\boldsymbol E_0$ and $\boldsymbol H_0$ are also assumed to be compactly
supported in $\Omega$. Due to the unbounded structure of the medium, it is no
longer valid to impose the usual Silver--M\"{u}ler radiation condition. We
employ the following radiation condition: the electromagnetic fields
$(\boldsymbol E, \boldsymbol H)$ consist of bounded outgoing waves in
$\Omega_j$. 

\subsection{Functional spaces}

We introduce some Sobolev space notation. For $ u \in L^2 (\Gamma_j)$, we denote
by $\hat u$ the  Fourier transform of $u$, i.e., 
 \[
 \hat u (\boldsymbol \xi, h_j)=\frac{1}{2 \pi} \int_{\mathbb R^2} u (\boldsymbol
\rho, h_j) e^{-{\rm i} \boldsymbol \rho \cdot \boldsymbol \xi} {\rm d}
\boldsymbol \rho,
 \]
where $\boldsymbol \xi= (\xi_1, \xi_2)^\top \in \mathbb R^2$ and $\boldsymbol
\rho =(x, y)^\top \in \mathbb R^2.$ Denote by $C_{\boldsymbol \rho}^{\infty}$
the linear space of infinitely differentiable functions with compact support
with respect to the variable $\boldsymbol \rho$ on $\Omega.$ Let $L^2
(\Omega)$ be the space of complex square integrable functions on  $\Omega$.
It follows from the Parseval identity that we have
\[
\|u\|^2_{L^2(\Omega)}=\int_{h_2}^{h_1} \int_{\mathbb R^2} |u(\boldsymbol \rho,
z)|^2 {\rm d} \boldsymbol \rho  {\rm  d} z=\int_{h_2}^{h_1} \int_{\mathbb
R^2}|\hat {u} (\boldsymbol \xi, z)|^2 {\rm d}\boldsymbol \xi  {\rm d} z.
\]

Introduce the functional spaces
\begin{align*}
\boldsymbol H({\rm curl}, \Omega)&=\{ \boldsymbol u \in \boldsymbol L^2(\Omega),
\nabla \times \boldsymbol u \in \boldsymbol L^2 (\Omega) \},\\
\boldsymbol H_0({\rm curl}, \Omega)&=\{ \boldsymbol u \in \boldsymbol
H({\rm curl}, \Omega), \boldsymbol{u}\times\boldsymbol n_j=0~\text{on} ~
\Gamma_j\},
\end{align*}
which are Sobolev spaces with the norm
\[
\|\boldsymbol u\|_{\boldsymbol H( {\rm curl}, \Omega)}=\Big( \| \boldsymbol
u\|^2_{\boldsymbol L^2 (\Omega)}+\|\nabla \times \boldsymbol u\|^2_{\boldsymbol
L^2 (\Omega)}\Big)^{1/2}.
\]
Given $\boldsymbol u = (u_1(\boldsymbol \rho, z), u_2 (\boldsymbol \rho, z), u_3
(\boldsymbol \rho, z))^\top\in\boldsymbol H ({\rm curl}, \Omega),$ it has
the inverse Fourier transform:
\[
\boldsymbol u(\boldsymbol \rho, z)= \frac{1}{2\pi}\int_{\mathbb R^2} \big(\hat
u_1(\boldsymbol \xi, z), \hat u_2 (\boldsymbol \xi, z), \hat u_3 (\boldsymbol
\xi, z)\big)^{\top} e^{{\rm i} \boldsymbol \rho \cdot \boldsymbol \xi} {\rm d}
\boldsymbol \xi.
\]
The norm in $\boldsymbol H ({\rm curl}, \Omega)$ can be defined via Fourier
coefficients:
\begin{align*}
\| \boldsymbol u\|^2_{\boldsymbol H ({\rm curl}, \Omega)}
=\int_{h_2}^{h_1} \int_{\mathbb R^2}
&\big[|\hat u_1 (\boldsymbol \xi, z)|^2 +|\hat u_2 (\boldsymbol \xi, z)|^2+|\hat
u_3 (\boldsymbol \xi, z)|^2+|{\rm i} \xi_{2} \hat u_3 (\boldsymbol \xi, z) -\hat
u_2'(\boldsymbol \xi, z)|^2\\
&+|\hat u_1'(\boldsymbol \xi, z)-{\rm i}\xi_1 \hat u_3 (\boldsymbol \xi,
z)|^2+|\xi_1 \hat u_2 (\boldsymbol \xi, z)- \xi_2 \hat u_1 (\boldsymbol \xi,
z)|^2 \big]{\rm d}\boldsymbol \xi {\rm d} z,
\end{align*}
where $\hat u_j'(\boldsymbol \xi, z)=\partial_z \hat u_j (\boldsymbol \xi, z).$

\begin{lemm}
$C_{\boldsymbol \rho}^{\infty} (\Omega)^3 $ is dense in $\boldsymbol H ({\rm
curl}, \Omega).$ 
\end{lemm}

\begin{proof}
Noting that $C_0^{\infty} (\mathbb R^3)^3$ is dense in $\boldsymbol H ({\rm
curl}, \mathbb R^3)$ , we have $C_0^{\infty} (\mathbb R^3)^3|_{\Omega}$ is dense
in $\boldsymbol H ({\rm curl}, \mathbb R^3)|_{\Omega}.$ From the Sobolev
extension theorem, $\boldsymbol H ({\rm curl}, \mathbb
R^3)|_{\Omega}=\boldsymbol H ({\rm curl}, \Omega).$ Therefore $C_{\boldsymbol
\rho}^{\infty} (\Omega)^3\supseteq C_0^{\infty} (\mathbb R^3)^3|_{\Omega}$ is
dense in  $\boldsymbol H ({\rm curl}, \Omega).$ 
\end{proof}

This density lemma is useful to deal with the infinite domain $\Omega$. We may
prove the results only on $ C_{\boldsymbol \rho}^{\infty} (\Omega)^3$ and then
extend them by limiting argument to more general functions such as those in
$\boldsymbol H (\rm curl, \Omega)$. Consequently, the boundary integrals only
on $\Gamma_j$ need to be considered when formulating the variational problems in
$\Omega$.

For any vector field $\boldsymbol u= (u_1, u_2, u_3)^{\top}$, denote by
\[
\boldsymbol u_{\Gamma_j}= \boldsymbol n_j \times( \boldsymbol u \times
\boldsymbol n_j)=(u_1 (x, y, h_j), u_2 (x, y, h_j), 0)^\top
\]
the tangential component on $\Gamma_j$, where $\boldsymbol n_1=(0, 0, 1)^\top$
and $\boldsymbol n_2=(0, 0, -1)^\top$ are the unit outward normal vectors on
$\Gamma_1$ and $\Gamma_2$, respectively. For any smooth vector
$\boldsymbol u=(u_1, u_2, u_3)^\top$ defined on $\Gamma_j$, let ${\rm div
}_{\Gamma_j} \boldsymbol u =\partial_x u_1+\partial_y u_2$ and ${\rm
curl}_{\Gamma_j} \boldsymbol u =\partial_x u_2 -\partial_y u_1$ be the surface
divergence and surface scalar curl of the field $\boldsymbol u$. For a smooth
scalar function $u$, denote by $\nabla_{\Gamma_j} u =(\partial_x u, \partial_y
u, 0)^{\top} $ the surface gradient on $\Gamma_j$. 

Let $H^{-1/2} (\Gamma_j)$ be the completion of $L^2(\Gamma_j)$ in the norm 
\begin{align*}
\|u\|_{H^{-1/2}(\Gamma_j)} =\Big(\int_{\mathbb R^2}(1+|\boldsymbol
\xi|^2)^{-1/2} |\hat u|^2 {\rm d} \boldsymbol \xi\Big)^{1/2}.
\end{align*}
Introduce two tangential functional spaces:
\begin{align*}
\boldsymbol H^{-1/2}({\rm curl}, \Gamma_j) &=\{ \boldsymbol u \in H^{-1/2}
(\Gamma_j)^3: u_3=0, ~{\rm curl }_{\Gamma_j} \boldsymbol u \in H^{-1/2}
(\Gamma_j)\},\\
\boldsymbol H^{-1/2}({\rm div}, \Gamma_j) &=\{ \boldsymbol u \in
H^{-1/2}(\Gamma_j)^3: u_3=0, ~{\rm div}_{\Gamma_j} \boldsymbol u \in H^{-1/2}
(\Gamma_j)\},
\end{align*}
which are equipped with the norms: 
\begin{align*}
\|\boldsymbol u\|_{\boldsymbol H^{-1/2}({\rm curl},
~\Gamma_j)}&=\Big(\int_{\mathbb R^2} (1+|\boldsymbol \xi|^2)\big[ |\hat u_1|^2
+|\hat u_2|^2 +|\xi_1 \hat u_2 -\xi_2 \hat u_1|^2\big] {\rm d} \boldsymbol
\xi\Big)^{1/2},\\
\| \boldsymbol u\|_{\boldsymbol H^{-1/2} ({\rm div},~
\Gamma_j)}&=\Big(\int_{\mathbb R^2} (1+|\boldsymbol \xi|^2)\big[ |\hat u_1|^2
+|\hat u_2|^2 +|\xi_1 \hat u_1 +\xi_2 \hat u_2|^2\big] {\rm d} \boldsymbol
\xi\Big)^{1/2}. 
\end{align*}

The following two lemmas are concerned with the duality between the spaces 
$\boldsymbol{H}^{-1/2}({\rm div}, \Gamma_j)$ and $ \boldsymbol
H^{1/2}({\rm curl}, \Gamma_j)$ and the trace regularity in $\boldsymbol H( {\rm
curl}, \Omega)$. The proofs can be found in \cite[Lemma 2.3, Lemma
2.4]{Li2011electromagnetic}. 

\begin{lemm}\label{AJ}
The spaces $\boldsymbol{H}^{-1/2}({\rm div}, \Gamma_j)$ and $ \boldsymbol
H^{1/2}({\rm curl}, \Gamma_j)$ are mutually adjoint with respect to
the scalar product in $L^2(\Gamma_j)^3$ defined by
\begin{equation}\label{ID}
\langle \boldsymbol u,  \boldsymbol v\rangle_{\Gamma_j}= \int_{\Gamma_j}
\boldsymbol u \cdot \bar {\boldsymbol v} {\rm d } \gamma_j=\int_{\mathbb R^2}
(\hat{u}_1\bar{\hat v}_1+\hat{u}_2\bar{\hat v}_2){\rm d}\boldsymbol\xi.
\end{equation}
\end{lemm}

\begin{lemm} \label{TaT}
Let $C=\max \{\sqrt {1+ (h_1-h_2)^{-1}}, \sqrt 2\}$. We have the estimate
\[
\|\boldsymbol u\|_{\boldsymbol H^{-1/2}({\rm curl}, ~\Gamma_j)} \leq C
\|\boldsymbol u\|_{\boldsymbol H ({\rm curl}, \Omega)},\quad\forall ~
\boldsymbol u \in \boldsymbol H( {\rm curl}, \Omega).
\]
\end{lemm}

Next we introduce some properties of the Laplace transform. Let $s=s_1+{\rm i}
s_2$ with $s_1>0, s_2\in\mathbb{R}$. Define by $\breve {\boldsymbol u} (s)$
the Laplace transform of $\boldsymbol u(t)$, i.e.,
\[
\breve {\boldsymbol u}(s) =\mathscr L (\boldsymbol u) (s)=\int_0^{\infty} e^{-s
t} \boldsymbol u (t) {\rm d} t.
\]
Using the integration by parts yields
\begin{equation}\label{A1}
\int_0^t \boldsymbol u (\tau) {\rm d} \tau =\mathscr
L^{-1} (s^{-1} \breve {\boldsymbol u} (s)),
\end{equation}
where $\mathscr L^{-1}$ is the inverse Laplace transform. It is also easy to
verify that  
\begin{equation}\label{a2}
\boldsymbol u(t)=\mathscr F^{-1} \big( e^{s_1 t} \mathscr L (\boldsymbol u) (s_1+s_2)\big),
\end{equation}
where $\mathscr F ^{-1}$ denotes the inverse Fourier transform with respect to
$s_2$. Recall the Plancherel or the Parseval identity for the Laplace transform
(cf. \cite[(2.46)]{Cohen2007}):
\begin{equation}\label{PI}
 \frac{1}{2 \pi} \int_{- \infty}^{\infty} \breve {\boldsymbol  u} (s) \breve
{\boldsymbol   v} (s) {\rm d} s_2= \int_0^{\infty} e^{- 2 s_1 t} {\boldsymbol
u} (t) {\boldsymbol v} (t) {\rm d} t,\quad\forall ~ s_1>\lambda,
\end{equation}
where $\breve  {\boldsymbol u}= \mathscr L (\boldsymbol u), \breve {\boldsymbol
v}= \mathscr L (\boldsymbol v)$ and $\lambda$ is abscissa of convergence for the
Laplace transform of $\boldsymbol u$ and $\boldsymbol v.$

The following lemma (\cite[Theorem 43.1]{Treves 1975}) is an analogue of the
Paley--Wiener--Schwarz theorem for the Fourier transform of distributions with
compact support in the case of the Laplace transform.

\begin{lemm}\label {A2}
 Let $\breve {\boldsymbol h} (s)$ be a holomorphic function in the half-plane
$s_1 > \sigma_0$ and be valued in the Banach space $\mathbb E$. The following
two statements are equivalent:
\begin{enumerate}

\item There is a distribution $ \breve {\boldsymbol h} \in \mathcal
D_{+}'(\mathbb E)$ whose Laplace transform is equal to $\breve{\boldsymbol
h}(s)$;

\item There is a real $\sigma_1$ with $\sigma_0 \leq \sigma_1 <\infty$ and an
integer $m \geq 0$ such that for all complex numbers $s$ with ${\rm Re} s =s_1 >
\sigma_1,$ it holds that $\| \breve {\boldsymbol  h} (s)\|_{\mathbb E} \lesssim (1+|s|)^{m}$,

\end{enumerate}
where $\mathcal D'_{+}(\mathbb E)$ is the space of distributions on the real
line which vanishes identically in the open negative half line.
\end{lemm}

\subsection{Transparent boundary condition}

We introduce an exact time-domain TBC to formulate the scattering problem into
the following initial-boundary value problem:
\begin{equation}\label{EHP}
\begin{cases}
\nabla \times \boldsymbol E +\mu \partial_t \boldsymbol H=0, \quad \nabla \times
\boldsymbol H- \varepsilon \partial_t \boldsymbol E= \boldsymbol J \quad
&{\rm in} ~\Omega, ~ t>0,\\
\boldsymbol E|_{t=0} =\boldsymbol E_0, \quad \boldsymbol H|_{t=0}
=\boldsymbol H_0 \quad &{\rm in}~\Omega,\\
\mathscr {T}_j [\boldsymbol{E}_{\Gamma_j}]=\boldsymbol H \times n_j \quad &{\rm
on}~ \Gamma_j, ~ t>0,
\end{cases}
\end{equation}
where $\boldsymbol{E}_{\Gamma_j}$ is the tangential component of
$\boldsymbol{E}$ on $\Gamma_j$ and $\mathscr T_j$ is the time-domain
electric-to-magnetic capacity operator.

In what follows, we shall derive the formulation of the operators
$\mathscr T_j$ and show some of their properties. Since the
derivation of $\mathscr T_1$ and $\mathscr T_2$ is analogous, we will only show
the details for $\mathscr T_1$ and state the corresponding result on $\mathscr
T_2$ without derivation. 

Notice that $\boldsymbol J$ is supported in $\Omega$ and
$\varepsilon=\varepsilon_1, \mu=\mu_1$ in $\Omega_1$, the system of Maxwell
equations \eqref{EH} reduce to 
\begin{equation}\label{SE}
\nabla \times \boldsymbol E+ \mu_1\partial_t \boldsymbol H=0, \quad  \nabla
\times \boldsymbol H -\varepsilon_1 \partial_t \boldsymbol E=0 \quad
{\rm in}~\Omega_1, ~ t>0.
\end{equation}
Let $\breve {\boldsymbol  E} (\boldsymbol x, s)$ and $\breve {\boldsymbol H}
(\boldsymbol x, s)$ be the Laplace transform of $\boldsymbol E (\boldsymbol x,
t)$ and $\boldsymbol H (\boldsymbol x, t)$. Recall that
\[
\mathscr L (\partial_t \boldsymbol E) = s \breve {\boldsymbol E}
-{\boldsymbol E}_0, \quad \mathscr L (\partial_t \boldsymbol H)
=s \breve {\boldsymbol H} -\boldsymbol H_0.
\]
Taking the Laplace transform of \eqref{SE}, and noting that $\boldsymbol E_0$
and $\boldsymbol H_0$ are supported in $\Omega,$ we obtain the Maxwell
equations in the $s$-domain:
\begin{equation}\label{CME}
\nabla \times \breve {\boldsymbol E} + \mu_1 s \breve {\boldsymbol H}=0, \quad 
\nabla \times \breve {\boldsymbol H  } -\varepsilon_1  s \breve {\boldsymbol
E}=0 \quad {\rm in}~ \Omega_1, ~ s_1>0, ~ s_2\in\mathbb{R}. 
\end{equation}

Let $\breve {\boldsymbol E} =(\breve E_1, \breve E_2, \breve E_3)^{\top}$ and
$\breve {\boldsymbol H} =(\breve H_1, \breve H_2, \breve H_3)^{\top}$. Denote
by $\breve {\boldsymbol E}_{\Gamma_1} = (\breve {E}_1 (\boldsymbol \rho, h_1),
\breve {E}_2 (\boldsymbol \rho, h_1), 0)^{\top}$ the tangential component of the
electric field on $\Gamma_1$. Let $\breve {\boldsymbol H} \times \boldsymbol
n_1=(\breve H_2 (\boldsymbol \rho, h_1), -\breve H_1 (\boldsymbol \rho, h_1),
0)^{\top}$ be the tangential trace of the magnetic field on $\Gamma_1.$ It
follows from \eqref{CME} that
\begin{align*}
\breve {H}_2 (\boldsymbol \rho, h_1)&=\frac{1}{\mu_1 s}[\partial_x
\breve{E}_3(\boldsymbol \rho, h_1)-\partial_z \breve {E}_1 (\boldsymbol \rho,
h_1)],\\
 -\breve {H}_1 (\boldsymbol \rho, h_1)&=\frac{1}{\mu_1 s}[\partial_y
\breve{E}_3(\boldsymbol \rho, h_1)-\partial_z \breve {E}_2 (\boldsymbol \rho,
h_1)].
\end{align*}
Taking the Fourier transform of the above equations with respect to $\boldsymbol
\rho$ gives 
\begin{subequations}\label{FT}
\begin{align}
\hat {\breve H}_2 (\boldsymbol \xi, h_1)&=\dfrac {1}{\mu_1 s} [{\rm i}\xi_1 \hat
{\breve E}_3 (\boldsymbol \xi, h_1)-\partial_z \hat {\breve E}_1 (\boldsymbol
\xi, h_1)],\\
-\hat {\breve H}_1 (\boldsymbol \xi, h_1)&=\dfrac {1}{\mu_1 s} [{\rm i}\xi_2
\hat {\breve E}_3 (\boldsymbol \xi, h_1)-\partial_z \hat {\breve E}_2
(\boldsymbol \xi, h_1)].
\end{align}
\end{subequations}

Observe that the medium is homogeneous in $\Omega_1$, which gives $\nabla \cdot
\breve {\boldsymbol E}=0$ in $\Omega_1.$ Eliminating the magnetic field from
\eqref{CME} and using the divergence free condition in $\Omega_1$, we obtain
the Helmholtz equation for the components of the electric field: 
\[
\begin{cases}
\Delta \breve {E}_j (\boldsymbol \rho,  z)-\varepsilon_1 \mu_1 s^2 \breve {E}_j
(\boldsymbol \rho, z)=0 \quad &{\rm in}~ \Omega_1,\\
\breve {E}_j (\boldsymbol \rho, z)= \breve E_j (\boldsymbol \rho, h_1) \quad
&{\rm on}~ \Gamma_1.
\end{cases}
\]
Taking the Fourier transform with respect to $\boldsymbol\rho$ of the above
equations yields
\[
\begin{cases}
\partial_z ^2 \hat {\breve E}_j  - (\varepsilon_1 \mu_1 s^2 +|\boldsymbol
\xi|^2) \hat {\breve E}_j=0, \quad & z>h_1,\\
\hat {\breve E}_j =\hat {\breve E}_j (\boldsymbol \xi, h_1), \quad & z=h_1.
\end{cases}
\]
Solving the above equations and using the bounded outgoing condition, we
obtain the solution:
\begin{equation}\label{ES}
\hat {\breve E}_j (\boldsymbol \xi, z)=\hat {\breve E}_j (\boldsymbol \xi, h_1)
e^{-\beta_1(\boldsymbol \xi) (z-h_1)}, \quad z>h_1,
\end{equation}
where
\[
\beta_1(\boldsymbol \xi)= (\varepsilon_1 \mu_1 s^2 +|\boldsymbol \xi|^2)^{1/2},
\quad {\rm Re} \beta_1 (\boldsymbol \xi)>0.
\]
Taking the derivative of \eqref{ES} with respect to $z$ and
evaluating it at $z=h_1$, we get 
\[
\partial_z \hat {\breve E}_j (\boldsymbol \xi, h_1)=-\beta_1 (\boldsymbol \xi)
\hat {\breve E}_j (\boldsymbol \xi, h_1).
\]
Noting that $\nabla \cdot \breve {\boldsymbol {E}}=\partial_x \breve {E}_1
+\partial_y \breve {E}_2 +\partial_z \breve E_3=0$ in $\Omega_1$ and
$\beta_1 (\boldsymbol \xi) \neq 0$ for all $\boldsymbol \xi,$ we deduce that
\[
\hat {\breve E}_3 (\boldsymbol \xi, h_1)= \frac{-1}{\beta_1 (\boldsymbol \xi)}
\partial_z \hat {\breve E}_3 (\boldsymbol \xi, h_1)=\frac{{\rm
i}}{\beta_1(\boldsymbol \xi)}[\xi_1 \hat {\breve E}_1(\boldsymbol \xi,
h_1)+\xi_2 \hat {\breve E}_2(\boldsymbol \xi, h_1)].
\]
Therefore, we have from \eqref{FT} that 
\begin{align*}
\hat {\breve H}_2(\boldsymbol \xi, h_1)&=\frac{1}{\mu_1 s}\left[ 
\frac{-\xi_1}{ \beta_1 (\boldsymbol \xi)} \bigl(\xi_1 \hat {\breve
E}_1(\boldsymbol \xi, h_1)+\xi_2 \hat {\breve E}_2(\boldsymbol \xi,
h_1)\bigr)+\beta_1(\boldsymbol \xi) \hat{ \breve {E}}_1(\boldsymbol \xi, h_1)
\right],\\
-\hat {\breve H}_1(\boldsymbol \xi, h_1)&=\frac{1}{\mu_1 s} \left[
\frac{-\xi_2}{\beta_1 (\boldsymbol \xi)} \bigl(\xi_1 \hat {\breve
E}_1(\boldsymbol \xi, h_1)+\xi_2 \hat {\breve E}_2 (\boldsymbol \xi, h_1)\bigr)
+\beta_1 (\boldsymbol \xi) \hat {\breve E}_2 (\boldsymbol \xi, h_1)\right],
\end{align*}
or equivalently,
\begin{align*}
\hat {\breve H}_2(\boldsymbol \xi, h_1)&=\frac{1}{\mu_1 s \beta_1(\boldsymbol
\xi)}\left[\varepsilon_1 \mu_1 s^2 \hat {\breve E}_1(\boldsymbol \xi,
h_1)+\xi_2\bigl(\xi_2 \hat {\breve E}_1(\boldsymbol \xi, h_1)-\xi_1 \hat {\breve
E}_2 (\boldsymbol \xi, h_1 )\bigr)\right],\\
-\hat {\breve H}_1(\boldsymbol \xi, h_1)&=\frac{1}{\mu_1 s \beta_1(\boldsymbol
\xi)} \left[\varepsilon_1 \mu_1 s^2 \hat {\breve E}_2 (\boldsymbol \xi, h_1)+
\xi_1 \bigl(\xi_1 \hat {\breve E}_2(\boldsymbol \xi, h_1) -\xi_2 \hat {\breve
E}_1(\boldsymbol \xi, h_1)\bigr)\right].
\end{align*}

For any tangential vector $\boldsymbol u= (u_1, u_2, 0)^{\top}$ on $\Gamma_1$, 
define the capacity operator $\mathscr B_1$:
\[
\mathscr B_1 [\boldsymbol u]=(v_1, v_2, 0)^{\top},
\]
where
\begin{subequations}\label{CO1}
\begin{align}
\hat v_1&=\frac{1}{\mu_1 s}\big[-\frac{\xi_1}{ \beta_1}(\xi_1 \hat u_1 +\xi_2
\hat u_2)+ \beta_1 \hat u_1\big],\\
\hat v_2&=\frac{1}{ \mu_1 s} \big[- \frac{\xi_2}{\beta_1}(\xi_1 \hat u_1+\xi_2
\hat u_2)+\beta_1 \hat u_2 \big],
\end{align}
\end{subequations}
or equivalently,
\begin{subequations}\label{CO2}
\begin{align}
\hat v_1&=\frac{1}{\mu_1 s \beta_1}\big[ \varepsilon_1 \mu_1 s^2 \hat u_1 +
\xi_2 (\xi_2 \hat u_1 -\xi_1 \hat u_2)\big],\\
\hat v_2&= \frac{1}{ \mu_1 s \beta_1 } \big [ \varepsilon_1 \mu_1 s^2 \hat u_2
+\xi_1 (\xi_1 \hat u_2 -\xi_2 \hat u_1) \big].
\end{align}
\end{subequations}
Similarly, for any tangential vector $\boldsymbol u= (u_1, u_2, 0)$ on
$\Gamma_2$, define the capacity operator $\mathscr B_2$:
\[
\mathscr B_2 [\boldsymbol u] =(v_1, v_2, 0)^{\top},
\]
where
\begin{subequations}\label{CO3}
\begin{align}
\hat v_1&=\frac{1}{\mu_2 s}\big[\beta_2 \hat u_1 -\frac{\xi_1}{ \beta_2}(\xi_1
\hat u_1+\xi_2 \hat u_2) \big],\\
\hat v_2&=\frac{1}{ \mu_2 s} \big[ \beta_2 \hat u_2 -\frac{\xi_2}{ \beta_2}
(\xi_1 \hat u_1+\xi_2 \hat u_2) \big],
\end{align}
\end{subequations}
or equivalently,
\begin{subequations}\label{CO4}
\begin{align}
\hat v_1&=\frac{1}{\mu_2 s \beta_2}\big[ \varepsilon_2 \mu_2 s^2 \hat u_1 +
\xi_2 (\xi_2 \hat u_1 -\xi_1 \hat u_2)\big],\\
\hat v_2&= \frac{1}{ \mu_2 s \beta_2 } \big [ \varepsilon_2 \mu_2 s^2 \hat u_2
+\xi_1 (\xi_1 \hat u_2 -\xi_2 \hat u_1) \big],
\end{align}
\end{subequations}
where
\[
\beta_2(\boldsymbol \xi)=(\varepsilon_2 \mu_2 s^2+|\boldsymbol \xi|^2)^{1/2},
\quad {\rm Re} \beta_2 (\boldsymbol \xi) >0.
\]

For any vector field $\breve {\boldsymbol E} \in \boldsymbol H ({\rm curl},
\Omega)$, it follows from Lemma \ref{TaT} that its tangential component $\breve
{\boldsymbol E}_{\Gamma_j} \in \boldsymbol H^{-1/2} ({\rm curl}, \Gamma_j).$
Using the capacity operators, we may propose the following TBC in the
$s$-domain:
\begin{equation}\label{TBC}
\mathscr B_j [\breve {\boldsymbol E}_{\Gamma_j}]=\breve {\boldsymbol H} \times
\boldsymbol n_j \quad {\rm on} ~ \Gamma_j,
\end{equation}
where the capacity operator $\mathscr B_j$ maps the tangential component of the
electric field to the tangential trace of  the magnetic field. Taking the
inverse Laplace transform of \eqref{TBC} yields the TBC in the time-domain:
\[
\mathscr T_j [\boldsymbol E_{\Gamma_j}]= \boldsymbol H \times  \boldsymbol
n_j,
\]
where $\mathscr T_j=\mathscr L^{-1} \circ \mathscr B_j \circ \mathscr
L$. Equivalently, we may eliminate the magnetic field and obtain an alternative
TBC for the electric field in the $s$-domain:
\begin{equation}\label{ATBC}
\mu_j ^{-1 } s^{-1 }(\nabla \times \breve {\boldsymbol E}) \times \boldsymbol 
n_j +\mathscr B_j [\breve {\boldsymbol E}_{\Gamma_j}]=0 \quad {\rm on
}~ \Gamma_j.
\end{equation}
Correspondingly, by taking the inverse Laplace transform of (\ref{ATBC}), we may
derive an alternative TBC for the electric field in the time-domain:
\begin{equation}\label{ATTBC}
\mu_j^{-1} (\nabla \times \boldsymbol E) \times \boldsymbol n_j +\mathscr  C_j
[\boldsymbol E_{\Gamma_j}]=0 \quad {\rm on}\quad \Gamma_j,
\end{equation}
where $\mathscr C_j = \mathscr L^{-1} \circ s \mathscr B_j \circ \mathscr L.$

\begin{lemm}\label{TC}
The capacity operator $\mathscr B_j :\boldsymbol H^{-1/2} ({\rm
curl}, \Gamma_j) \rightarrow \boldsymbol H^{-1/2} ({\rm div}, \Gamma_j)$ is
continuous.
\end{lemm}

\begin{proof}
For any $\boldsymbol u =(u_1, u_2 ,0)^{\top}, \boldsymbol w =(w_1, w_2,
0)^{\top} \in \boldsymbol H^{-1/2} ({\rm curl}, \Omega),$ let $\mathscr B_j
\boldsymbol u = ( v_1, v_2, 0)^{\top}.$ It follows from the definitions
\eqref{ID}, \eqref{CO2}, and \eqref{CO4} that
\begin{align*}
\langle \mathscr B_j  \boldsymbol u, \boldsymbol w \rangle_{\Gamma_j}
&= \int_{\mathbb R^2} (\hat v_1 \bar{\hat w}_1 +\hat v_2 \bar {\hat w}_2) {\rm
d} \boldsymbol \xi \\
&= \int_{\mathbb R^2} \frac{1}{ \mu_j s \beta_j} \big[ \varepsilon_j \mu_j s^2
(\hat u_1 \bar {\hat w}_1 +\hat u_2 \bar {\hat w}_2) +(\xi_1 \hat u_2 -\xi_2
\hat u_1) (\xi_1 \bar {\hat w}_2- \xi_2 \bar {\hat w}_1)\big]{\rm d} \boldsymbol
\xi\\
&=\int_{\mathbb R^2} \frac{(1+|\boldsymbol \xi|^2)^{1/2}}{\mu_j s \beta_j}
(1+|\boldsymbol \xi|^2)^{-1/2}\big[\varepsilon_j \mu_j s^2 (\hat u_1 \bar {\hat
w}_1 +\hat u_2 \bar {\hat w}_2) +(\xi_1 \hat u_2 -\xi_2 \hat u_1) (\xi_1 \bar
{\hat w}_2- \xi_2 \bar {\hat w}_1)  \big] {\rm d} \boldsymbol \xi.
\end{align*}
To prove the lemma, it is required to estimate
\[
 \frac{(1+|\boldsymbol \xi|^2)^{1/2}}{|\beta_j|}.
\]
Let
\[
\varepsilon_j \mu_j s^2 =a_j+ {\rm i} b_j,
\]
where 
\[
a_j=\varepsilon_j \mu_j (s_1^2-s_2^2),\quad b_j =2 \varepsilon_j \mu_j s_1 s_2. 
\]
Denote
\[
\beta_j^2=\varepsilon_j\mu_j s^2+|\boldsymbol \xi|^2=\phi_j +{\rm i} b_j,
\]
where 
\[
\phi_j ={\rm Re} (\varepsilon_j \mu_j s^2)+|\boldsymbol \xi|^2=a_j
+|\boldsymbol \xi|^2.
\]
A simple calculation gives
\[
\frac{(1+|\boldsymbol \xi|^2)^{1/2}}{|\beta_j|}=\left[  \frac{(1+\phi_j -a_j)^2}{\phi_j^2+b_j^2}\right]^{1/4}.
\]
Let
\[
F_j (t)=\frac{(1+t -a_j)^2}{t^2+b_j^2},
\]
which gives 
\[
F_j'(t)= \frac{2 (1+t -a_j) (b_j^2 -t(1-a_j))}{(t^2+b_j^2)^2}.
\]
We consider three cases:
\begin{itemize}

\item[(i)] $1-a_j> 0$. It can be verified that the function $F_j(t)$ increases
for $ a_j\leq t \leq K_j=b_j^2/(1-a_j)$ and decreases for  $t> K_j$. Hence
$F_j(t)$ reaches its maximum at $t=K_j$, i.e., 
\[
\frac{(1+\phi_j -a_j)^2}{\phi_j^2+b_j^2}=F_j(\phi_j) \leq
F_j(K_j)=\frac{(1-a_j)^2+b_j^2}{b_j^2}.
\]

\item[(ii)] $1-a_j=0$. It is easy to verify 
\[
F_j (t)=\frac{t^2}{t^2+b_j^2} \leq 1.
\]
which yields that 
\[
F_j(\phi_j) \leq 1 \leq \frac{(1-a_j)^2+b_j^2}{b_j^2}.
\]

\item[(iii)] $1-a_j <0$. It follows from $K_j \leq a_j$ that $F_j(t)$ increases
for $t \leq K_j$ and decreases for  $ K_j < t$. Since $\phi_j =a_j
+|\boldsymbol \xi|^2 \geq a_j$, we have 
\[
F_j(\phi_j) \leq F_j(a_j)=\frac{1}{a_j^2 +b_j^2} \leq
F_j(K_j)=\frac{(1-a_j)^2+b_j^2}{b_j^2}.
\]
\end{itemize}
Combing the above estimates yields
\[
|\langle \mathscr B_j \boldsymbol u,  \boldsymbol w\rangle_{\Gamma_j}| \leq C_j
\|\boldsymbol u\|_{\boldsymbol H^{-1/2} ({\rm curl},~ \Gamma_j)} \|\boldsymbol
w\|_{\boldsymbol H^{-1/2} ({\rm curl},~ \Gamma_j)}
\]
where
\[
C_j= \frac{1}{\mu_j s_1}\left[\frac{(1-a_j)^2+b_j^2}{b_j^2}\right]^{1/4}\times
\max \{(a_j^2 +b_j^2)^{1/2}, ~ 1\}.
\]
Following from Lemma \ref{AJ}, we have
\[
\|\mathscr B_j \boldsymbol u\|_{\boldsymbol H^{-1/2} ({\rm div},~ \Gamma_j)}
\leq C \sup \limits_{\boldsymbol w \in \boldsymbol H^{-1/2}({\rm curl},~
\Gamma_j)}\frac{|\langle \mathscr B_j \boldsymbol u,  \boldsymbol
w\rangle_{\Gamma_j}|}{\|\boldsymbol w\|_{\boldsymbol H^{-1/2}({\rm curl},~
\Gamma_j)}} \leq C C_j \|\boldsymbol u\|_{\boldsymbol H^{-1/2}({\rm curl},~
\Gamma_j)},
\]
which completes the proof. 
\end{proof}

\begin{lemm}\label{TP}
We have 
\[
{\rm Re} \langle \mathscr B_j \boldsymbol u, \boldsymbol u \rangle_{\Gamma_j}
\geq 0, \quad \forall ~\boldsymbol u \in \boldsymbol H^{-1/2 }({\rm curl},
\Gamma_j).
\]
\end{lemm}

\begin{proof}
By definitions (\ref{ID}), (\ref{CO1}), and (\ref{CO3}), we obtain
\begin{align*}
\langle \mathscr B_j \boldsymbol u, \boldsymbol u \rangle_{\Gamma_j}
&=\frac{1}{ \mu_j s} \int_{\mathbb R^2} \big[ \beta_j (|\hat u_1|^2 +|\hat u_2|^2) -
\frac{1}{\beta_j}
|\xi_1 \hat u_1 +\xi_2 \hat u_2|^2 \big]{\rm d} \boldsymbol
\xi\\
&=\frac{1}{\mu_j |s|^2}\int_{\mathbb R^2} \big[ \bar s \beta_j (|\hat
u_1|^2+|\hat u_2|^2)-\frac{ \bar s \bar {\beta_j}}{|\beta_j|^2}|\xi_1 \hat u_1
+\xi_2 \hat u_2|^2 \big] {\rm d} \boldsymbol \xi.
\end{align*}
Let $\beta_j =m_j +{\rm i}  n_j$ with $m_j >0.$  Taking the real part of the
above equation gives
\begin{align*}
{\rm Re} \langle \mathscr B_j \boldsymbol u, \boldsymbol u \rangle_{\Gamma_j}  =
\frac{1}{ \mu_j |s|^2 } \int_{\mathbb R^2} \big[ &(m_j s_1 +n_j s_2) (|\hat
u_1|^2+|\hat u_2|^2)\\
-&\frac{(m_j s_1- n_j s_2) }{m_j^2+n_j^2} |\xi_1 \hat u_1
+\xi_2 \hat u_2|^2\big] {\rm d} \boldsymbol \xi.
\end{align*}
Recalling  $\beta_j^2 =\varepsilon_j \mu_j s^2+|\boldsymbol \xi|^2,$  we have
\begin{align}
m_j^2 -n_j^2&= \varepsilon_j \mu_j (s_1^2 -s_2^2)+|\boldsymbol \xi|^2, \label
{R2}\\
 m_j n_j&= \varepsilon_j \mu_j s_1 s_2. \label{R3}
\end{align}
Using \eqref {R3}, we get
\begin{align*}
m_j s_1 +n_j s_2=\frac{s_1}{m_j} [m_j^2+\varepsilon_j \mu_j s_2^2],\quad
m_j s_1 -n_j s_2= \frac{s_1}{m_j} [m_j^2 -\varepsilon_j \mu_j s_2^2].
\end{align*}
If $ m_j^2 -\varepsilon_j \mu_j s_2^2 \leq 0$, we obtain
\begin{align*}
{\rm Re} \langle \mathscr B_j \boldsymbol u, \boldsymbol u \rangle_{\Gamma_j}  =
\frac{1}{ \mu_j |s|^2 } \int_{\mathbb R^2} &\frac{s_1}{m_j}\big[
(m_j^2+\varepsilon_j \mu_j s_2^2) (|\hat u_1|^2+|\hat u_2|^2) \\
-&\frac{(m_j^2-\varepsilon_j \mu_j s_2^2)}{m_j^2+n_j^2} |\xi_1 \hat u_1 +\xi_2
\hat u_2|^2\big] {\rm d} \boldsymbol \xi \geq 0.
\end{align*}
If $ m_j^2 -\varepsilon_j \mu_j s_2^2  >0$,  we have from the
Cauchy--Schwarz inequality that 
\[
\frac{(m_j^2 -\varepsilon_j \mu_j s_2^2)}{m_j^2+n_j^2} |\xi_1 \hat u_1 +\xi_2
\hat u_2|^2   \leq  \frac{ (m_j^2 -\varepsilon_j \mu_j s_2^2)}{ m_j^2+n_j^2}
|\boldsymbol\xi|^2(|\hat u_1|^2 +|\hat u_2 |^2),
\]
which gives
\begin{align}\label{R5}
{\rm Re} \langle \mathscr B_j \boldsymbol u, \boldsymbol u \rangle_{\Gamma_j}
\geq \frac{1}{\mu_j |s|^2} \int_{\mathbb R^2}& \frac{s_1}{m_j}\Bigl[ (m_j^2
+\varepsilon_j \mu_j s_2^2)\notag\\
-&\frac{(m_j^2 -\varepsilon_j \mu_j s_2^2)}{m_j^2+n_j^2}|\boldsymbol\xi|^2\Bigr]
(|\hat u_1|^2 +|\hat u_2|^2) {\rm d } \boldsymbol \xi.
\end{align}
Substituting \eqref {R2} into \eqref{R5} yields
\begin{align*}
{\rm Re} \langle \mathscr B_j \boldsymbol u, \boldsymbol u \rangle_{\Gamma_j}
\geq \frac{1}{\mu_j |s|^2} \int_{\mathbb R^2}&\frac{s_1}{m_j(m_j^2+n_j^2)}
\bigl[(m_j ^2 +\varepsilon_j \mu_j s_2^2)(n_j^2+\varepsilon_j \mu_j s_2^2)\\
 +&(m_j^2 -\varepsilon_j \mu_j s_2^2 )(n_j^2+\varepsilon_j \mu_j s_1^2)\bigr]
(|\hat u_1|^2+|\hat u_2|^2)  {\rm d} \boldsymbol \xi \geq 0,
\end{align*}
which completes the proof. 
\end{proof}

In the forthcoming sections, we shall use the method of energy to prove the
well-posedness and stability of the reduced problem (\ref{EHP}). We point out
that the method has also been adopted in \cite{LiWangWood2015} for solving the
time-dependent electromagnetic scattering problem from an open cavity.

\section{Two auxiliary problems}\label{at}

In this section, we present the energy estimates for two auxiliary problems, one
is the time-harmonic Maxwell equations with a complex wavenumber and another is
the time-domain Maxwell equations with a perfectly electrically conducting (PEC)
boundary condition. These estimates will be used for the proof of the main
results for the time-domain Maxwell equations \eqref{EHP}.

\subsection{Time-harmonic Maxwell's equations with a complex wavenumber}

We shall study the variational formulation for a time-harmonic Maxwell equations
with a complex wavenumber, which is a frequency version of the initial-boundary
value problem of the Maxwell equations under the Laplace transform. 

Consider the auxiliary boundary value problem:
\begin{equation}\label{AP}
\left \{
\begin{array}{ll}
\nabla \times \big((s \mu)^{-1} \nabla \times \boldsymbol u\big)+s \varepsilon
\boldsymbol u =\boldsymbol j \quad &\text{in}~ \Omega,\\
\mu_j ^{-1 } s^{-1 }(\nabla \times  {\boldsymbol u}) \times \boldsymbol  n_j
+\mathscr B_j [{\boldsymbol u}_{\Gamma_j}]=0\quad& \text{on}~ \Gamma_j,
\end{array}
\right.
\end{equation}
where $s= s_1 +{\rm i} s_2$ with $s_1>0, s_2 \in \mathbb R$ and $\boldsymbol{j}$
is assumed to be compactly supported in $\Omega.$

Multiplying the complex conjugate of a test function $\boldsymbol v \in
\boldsymbol H ({\rm curl}, \Omega),$ integrating over $\Omega$, and using
integration by parts, we arrive at the variational formulation  of \eqref{AP}:
Find $\boldsymbol u \in \boldsymbol H ({\rm curl}, \Omega)$ such that
\begin{equation}\label{ABV}
a_{\rm TH} (\boldsymbol u, \boldsymbol v) =\int_{\Omega} \boldsymbol j \cdot
\bar {\boldsymbol v} {\rm d } \boldsymbol x,  \quad \forall ~ \boldsymbol v \in
\boldsymbol H({\rm curl}, \Omega),
\end{equation}
where the sesquilinear form
\begin{equation}\label{ASF}
a_{\rm TH} (\boldsymbol u, \boldsymbol v) =\int_{\Omega} (s \mu)^{-1} (\nabla
\times \boldsymbol u) \cdot (\nabla \times \bar {\boldsymbol v}) {\rm d}
\boldsymbol x
+\int_{\Omega} s \varepsilon \boldsymbol u \cdot \bar {\boldsymbol v} {\rm d}
\boldsymbol x +\sum \limits_{j=1}^2 \langle \mathscr B_j [\boldsymbol
u_{\Gamma_j}], \boldsymbol v_{\Gamma_j} \rangle_{\Gamma_j}.
\end{equation}

\begin{theo}\label{AT}
The variational problem (\ref{ABV}) has a unique solution $\boldsymbol u \in
\boldsymbol H ( {\rm curl}, \Omega)$ which satisfies
\[
\|\nabla \times \boldsymbol u\|_{\boldsymbol L^2 (\Omega)}+ \|s \boldsymbol
u\|_{\boldsymbol L^2 (\Omega)} \lesssim s_1 ^{-1} \|s
\boldsymbol j\|_{\boldsymbol L^2
(\Omega)}.
\]
\end{theo}
\begin{proof}
 It suffices to show the coercivity of the sesquilinear form of $a_{\rm TH}$
since  the continuity follows directly from the Cauchy--Schwarz inequality,
Lemma \ref{TC}, and Lemma \ref{TaT}. 

Letting $\boldsymbol v =\boldsymbol u$, we have from \eqref{ASF} that 
\begin{equation}\label{AC}
a_{\rm TH} (\boldsymbol u, \boldsymbol u)= \int_{\Omega} (s \mu)^{-1} |\nabla
\times \boldsymbol u|^2 {\rm d} \boldsymbol x +\int_{\Omega} s \varepsilon
|\boldsymbol u|^2 {\rm d} \boldsymbol x +\sum \limits_{j=1}^2 \langle \mathscr
B_j [\boldsymbol u_{\Gamma_j}], \boldsymbol u_{\Gamma_j} \rangle_{\Gamma_j}.
\end{equation}
Taking the real part of \eqref{AC} and using Lemma \ref{TP}, we get
\begin{equation}\label{RP}
{\rm Re } a_{\rm TH} (\boldsymbol u, \boldsymbol u)\geq \frac{s_1}{|s|^2}
\big(\|\nabla \times \boldsymbol u\|^2_{\boldsymbol L^2 (\Omega)}
+\|s \boldsymbol u\|^2_{\boldsymbol L^2 (\Omega)}\big).
\end{equation}
It follows from the Lax--Milgram lemma that the variational problem \eqref{ABV} 
has a unique solution $\boldsymbol u  \in \boldsymbol H ({\rm curl}, \Omega)$.
Moreover, we have from \eqref{ABV} that
\begin{equation}\label{SE1}
|a_{\rm TH} (\boldsymbol u, \boldsymbol u)|\leq |s|^{-1}\|\boldsymbol
j\|_{\boldsymbol L^2 (\Omega)} \| s \boldsymbol u\|_{\boldsymbol L^2 (\Omega)}.
\end{equation}
Combing \eqref{RP}--\eqref{SE1} leads to
\[
\|\nabla \times \boldsymbol u\|^2_{\boldsymbol L^2 (\Omega)}+\|s \boldsymbol
u\|^2_{\boldsymbol L^2 (\Omega)} \lesssim s_1 ^{-1} \|s \boldsymbol 
j\|_{\boldsymbol L^2 (\Omega)} \| s \boldsymbol u\|_{\boldsymbol L^2 (\Omega)},
\]
which completes the proof after applying the Cauchy--Schwarz inequality.
\end{proof}

\subsection{Time-domain Maxwell's equations with PEC condition}

Consider the initial-boundary value problem for the time-domain Maxwell
equations with the PEC boundary condition: 
\begin{equation}\label{TME}
\begin{cases}
\nabla \times \boldsymbol U+ \mu \partial_t \boldsymbol V=0, \quad \nabla \times
\boldsymbol V- \varepsilon \partial_t \boldsymbol U=0 \quad &\text{in} ~ \Omega,
~ t>0,\\
\boldsymbol U \times \boldsymbol n_j=0 \quad  &\text{on} ~ \Gamma_j, ~
t>0,\\
\boldsymbol U|_{t=0}= \boldsymbol E_0, \quad \boldsymbol V|_{t=0}=\boldsymbol
H_0 \quad & \text{in} ~ \Omega,
\end{cases}
\end{equation}
where $\boldsymbol E_0, \boldsymbol H_0$ are assumed to be compactly supported 
in $\Omega$. 

Let $\breve {\boldsymbol U} =\mathscr L (\boldsymbol U)$ and $\breve
{\boldsymbol V} =\mathscr L (\boldsymbol V).$ Taking the Laplace transform of
\eqref{TME} and eliminating  $\breve {\boldsymbol V}$, we obtain the boundary
value problem:
\begin{equation}\label{BP}
\begin{cases}
\nabla \times \big( (s \mu)^{-1} \nabla \times \breve {\boldsymbol U} \big) + s
\varepsilon \breve {\boldsymbol U} = \breve {\boldsymbol j}\quad &\text{in} ~
\Omega,\\
\breve {\boldsymbol U} \times \boldsymbol n_j =0 \quad &\text{on}~\Gamma_j,
\end{cases}
\end{equation}
where $\breve {\boldsymbol j} =\varepsilon\boldsymbol E_0+s^{-1} \nabla \times
\boldsymbol H_0$. The variational formulation for \eqref {BP} is to find $\breve
{\boldsymbol U} \in \boldsymbol H_0 ({\rm curl, \Omega})$ such that
 \begin{equation}\label{APB}
 a_{\rm TD} (\breve {\boldsymbol U}, \boldsymbol v) =\int_{\Omega} \breve
{\boldsymbol j} \cdot \bar {\boldsymbol v} {\rm d} \boldsymbol x,
 \quad \forall~\boldsymbol v \in \boldsymbol H_0({\rm curl}, \Omega),
 \end{equation}
 where the sesquilinear from
 \[
 a_{\rm TD} ( \breve {\boldsymbol U}, \boldsymbol v ) =\int_{\Omega} (s
\mu)^{-1} (\nabla \times \breve {\boldsymbol U}) \cdot (\nabla \times \bar
{\boldsymbol v}) {\rm d} \boldsymbol x +\int_{\Omega} s \varepsilon \breve
{\boldsymbol U} \cdot \bar {\boldsymbol v} {\rm d} \boldsymbol x.
 \]
Following the same proof  as that in Theorem \ref{AT}, we may obtain the
well-posedness of the variation problem \eqref{APB} and its stability estimate. 

\begin{lemm}\label{ase}
The variational problem (\ref{BP}) has a unique solution $\breve {\boldsymbol U}
\in \boldsymbol H_0 ({\rm curl}, \Omega)$ which satisfies
\[
\|\nabla \times \breve {\boldsymbol U}\|_{\boldsymbol L^2 (\Omega)} +\|s  \breve
{\boldsymbol U}\|_{\boldsymbol L^2 (\Omega)}\lesssim  s_1 ^{-1} \| s \boldsymbol
E_0\|_{\boldsymbol L^2 (\Omega)} + s_1^{-1} \|\nabla \times \boldsymbol
H_0\|_{\boldsymbol L^2 (\Omega)}.
\]
\end{lemm}

 \begin{theo}
The auxiliary problem \eqref {TME} has  a unique solution $(\boldsymbol U,
\boldsymbol V)$, which satisfies the stability estimates:
\begin{align*}
 \|\boldsymbol U\|_{\boldsymbol L^2 (\Omega)}+\| \boldsymbol V \|_{\boldsymbol L
^2 (\Omega)}  &\lesssim \|\boldsymbol E_0\|_{\boldsymbol L^2 (\Omega)}
+\|\boldsymbol H_0\|_{\boldsymbol L ^2 (\Omega)},\\
  \| \partial_t \boldsymbol U\|_{\boldsymbol L^2 (\Omega)} +\| \partial_t
\boldsymbol V \|_{\boldsymbol L ^2 (\Omega)}   &\lesssim \|\nabla \times
\boldsymbol E_0\|_{\boldsymbol L^2 (\Omega)} +\|\nabla \times \boldsymbol
H_0\|_{\boldsymbol L ^2 (\Omega)},\\
  \| \partial^2_t \boldsymbol U\|_{\boldsymbol L^2 (\Omega)}+\| \partial^2_t
\boldsymbol V \|_{\boldsymbol L ^2 (\Omega)}  &\lesssim \|\nabla \times (\nabla
\times \boldsymbol E_0)\|_{\boldsymbol L^2 (\Omega)} +\|\nabla \times (\nabla
\times \boldsymbol H_0)\|_{\boldsymbol L ^2 (\Omega)}.
\end{align*}
\end{theo}
 
 \begin{proof}
Let $\breve {\boldsymbol U} =\mathscr L (\boldsymbol U)$ and $\breve
{\boldsymbol V}=\mathscr L (\boldsymbol V)$. Taking the Laplace transform of
\eqref{TME} and using the initial condition lead to
 \begin{equation}\label{LT}
 \left\{
 \begin{array}{ll}
 \nabla \times  \breve {\boldsymbol U} + s \mu \breve {\boldsymbol V} =\mu
\boldsymbol H_0, \quad \nabla \times \breve {\boldsymbol V} - s \varepsilon
\breve {\boldsymbol U}= - \varepsilon \boldsymbol E_0 \quad &\text{in}
~\Omega,\\
\breve {\boldsymbol U} \times \boldsymbol n_j =0 \quad &\text{on} ~ \Gamma_j.
 \end{array}
 \right.
 \end{equation}
 It follows from Lemma \ref{ase} that
 \[
\|\nabla \times \breve {\boldsymbol U}\|_{\boldsymbol L^2 (\Omega)} +\|s  \breve
{\boldsymbol U}\|_{\boldsymbol L^2 (\Omega)}\lesssim  s_1 ^{-1} \| s \boldsymbol
E_0\|_{\boldsymbol L^2 (\Omega)} + s_1^{-1} \|\nabla \times \boldsymbol
H_0\|_{\boldsymbol L^2 (\Omega)}. 
\]
Combing the above inequality and \eqref{LT} gives
\[
\|-s \mu \breve {\boldsymbol V} +\mu \boldsymbol H_0\|_{\boldsymbol L^2
(\Omega)} +\| \varepsilon ^{-1} \nabla \times \breve {\boldsymbol V} +
\boldsymbol E_0\| \lesssim  s_1 ^{-1} \| s \boldsymbol E_0\|_{\boldsymbol L^2
(\Omega)} + s_1^{-1} \|\nabla \times \boldsymbol H_0\|_{\boldsymbol L^2
(\Omega)}, 
\]
which shows that 
\[
\|\nabla \times {\breve {\boldsymbol V}}\|_{L^2 (\Omega)} +\| s \breve
{\boldsymbol V}\|_{\boldsymbol L^2 (\Omega)} \lesssim (1+s_1 ^{-1} |s|)
\|\boldsymbol E_0\|_{\boldsymbol L^2 (\Omega)} + \|\boldsymbol
H_0\|_{\boldsymbol L^2 (\Omega)} +s_1^{-1} \|\nabla \times \boldsymbol
H_0\|_{\boldsymbol L^2 (\Omega)}.
\]
It follows from \cite[Lemma 44.1]{Treves1975} that $\breve { \boldsymbol U}$ and
$\breve {\boldsymbol V}$ are holomorphic functions of $s$ on the half plane $s_1
>\bar\gamma>0,$  where $\bar \gamma$ is any positive constant. Hence we have
from Lemma \ref{A2} that the inverse Laplace transform of $\breve { \boldsymbol
U}$ and $\breve {\boldsymbol V}$ exist and they are  supported in $[0, \infty).$

Next we prove the stability by the energy function method. Define the energy  function
\[
e_1 (t)=\| \varepsilon^{1/2}  \boldsymbol U (\cdot, t)\|^2_{\boldsymbol
L^2(\Omega)} +\|\mu ^{1/2} \boldsymbol V (\cdot, t)\|^2_{\boldsymbol L^2
(\Omega)}.
\]
Using \eqref {TME} and integration by parts, we obtain
\begin{align*}
e_1 (t) -e_1 (0)
&=\int_{0}^{t} e_1'(\tau) {\rm d} \tau= 2 {\rm Re} \int_0^t \int_{\Omega} \big(
\varepsilon \partial_t \boldsymbol U \cdot \bar {\boldsymbol U}+ \mu \partial_t
\boldsymbol V \cdot \bar {\boldsymbol V}\big ) {\rm d} \boldsymbol x {\rm d}
\tau\\
&= 2{\rm Re} \int_0^t \int_{\Omega} (\nabla \times \boldsymbol V) \cdot \bar
{\boldsymbol U} -(\nabla \times {\boldsymbol U}) \cdot \bar {\boldsymbol V} {\rm
d} \boldsymbol x {\rm d} \tau\\
&= 2{\rm Re } \int_0^ t \int_{\Omega} \big[(\nabla \times  {\boldsymbol V})
\cdot \bar {\boldsymbol U} -(\nabla \times \bar {\boldsymbol V}) \cdot 
{\boldsymbol U} \big] {\rm d} \boldsymbol x  {\rm d} \tau -2 {\rm Re} \int_0^{t}
 \sum \limits_{j=1}^2 \langle \boldsymbol U \times \boldsymbol n_j, {\boldsymbol
V} \rangle_{\Gamma_j} {\rm d} \tau \\
&=0.
\end{align*}
Hence we have
\[
\| \varepsilon ^{1/2} \boldsymbol U (\cdot, t)\|^2_{\boldsymbol L^2} +\|\mu
^{1/2} \boldsymbol V (\cdot, t)\|^2_{\boldsymbol L^2 (\Omega)}=
\| \varepsilon ^{1/2} \boldsymbol E_0\|^2_{\boldsymbol L^2} +\|\mu ^{1/2}
\boldsymbol H_0\|^2_{\boldsymbol L^2 (\Omega)},
\]
which implies
\[
\|\boldsymbol U\|_{\boldsymbol L^2 (\Omega)} +\| \boldsymbol V \|_{\boldsymbol L ^2 (\Omega)}
  \lesssim \|\boldsymbol E_0\|_{L^2 (\Omega)} +\|\boldsymbol H_0\|_{\boldsymbol L ^2 (\Omega)}.
\]

Taking the first and second partial derivative of (\ref{TME}) with respect to $t$ yields
\begin{equation*}
\left\{
\begin{array}{ll}
\nabla \times \partial_t \boldsymbol U+ \mu \partial_t^2 \boldsymbol V=0, \quad
\nabla \times \partial_t \boldsymbol V -\varepsilon \partial_t^2 \boldsymbol U=0
\quad &\text{in} ~  \Omega, ~  t >0,\\
\partial_t \boldsymbol U \times \boldsymbol n_j =0 \quad  &\text{on} ~\Gamma_j,
~ t >0, \\
\partial_t \boldsymbol U|_{t=0}= \varepsilon^{-1} (\nabla \times \boldsymbol
H_0), \quad  \partial_t \boldsymbol V|_{t=0}= -\mu^{-1} \nabla \times
\boldsymbol E_0 \quad &\text{in}~\Omega
\end{array}
\right.
\end{equation*}
 and
 \begin{equation*}
\left\{
\begin{array}{ll}
\nabla \times \partial^2_t \boldsymbol U+ \mu \partial_t^3 \boldsymbol V=0,
\quad \nabla \times \partial^2_t \boldsymbol V -\varepsilon \partial_t^3
\boldsymbol U=0
\quad &\text{in} ~ \Omega, ~  t >0,\\
\partial^2_t \boldsymbol U \times \boldsymbol n_j =0 \quad  &\text{on}
~\Gamma_j, ~ t >0, \\
\partial^2_t \boldsymbol U|_{t=0}= - (\varepsilon \mu) ^{-1} (\nabla \times
(\nabla \times \boldsymbol E_0)), \quad&\text{in}~\Omega, \\
\partial^2_t \boldsymbol V|_{t=0}=-(\varepsilon \mu)^{-1} (\nabla \times
(\nabla \times \boldsymbol H_0)) \quad&\text{in}~\Omega.
\end{array}
\right.
\end{equation*}
Consider the energy functions
\[
e_2 (t)=\| \varepsilon^{1/2} \partial_t \boldsymbol U (\cdot,
t)\|^2_{\boldsymbol L^2 (\Omega)} +\|\mu ^{1/2} \partial_t \boldsymbol V (\cdot,
t)\|^2_{L^2(\Omega)}
\]
and
\[
e_3 (t)=\| \varepsilon^{1/2} \partial^2_t \boldsymbol U (\cdot,
t)\|^2_{\boldsymbol L^2 (\Omega)} +\|\mu ^{1/2} \partial^2_t \boldsymbol V
(\cdot, t)\|^2_{L^2(\Omega)}
\]
for the above two problems, respectively. Using the same steps for the first
inequality, we can derive the other two inequalities. The details are omitted. 
\end{proof}

\section{The reduced problem}\label{rp}

In this section, we present the main results of this work, which include the
well-posedness, stability, and a priori estimates for the scattering problem
\eqref{TME}. 

\subsection{Well-posedness}

Let $\boldsymbol e= \boldsymbol E -\boldsymbol U$ and $\boldsymbol h= 
\boldsymbol H- \boldsymbol V.$ Noting $\boldsymbol U\times \boldsymbol n_j=0$,
we have $\boldsymbol U_{\Gamma_j}=0$ and
$\mathscr{T}_j[\boldsymbol {U}_{\Gamma_j}]=0$. It follows from \eqref{EHP} and
\eqref{TME} that $\boldsymbol e$ and $\boldsymbol h$ satisfy the following
initial-boundary value problem:
\begin{equation}\label{nrp}
\begin{cases}
\nabla \times \boldsymbol e +\mu \partial_t \boldsymbol h =0, \quad \nabla
\times \boldsymbol h- \varepsilon \partial_t \boldsymbol e =\boldsymbol J \quad
 &{\rm in}~\Omega, ~ t>0,\\
\boldsymbol e|_{t=0} =0, \quad \boldsymbol h|_{t=0}=0 \quad
&{\rm in}~\Omega,\\
\mathscr T_j [\boldsymbol e_{\Gamma_j}]= \boldsymbol h \times \boldsymbol n_j
+\boldsymbol V \times \boldsymbol n_j \quad&{\rm on}~\Gamma_j, ~ t>0. 
\end{cases}
\end{equation}
Let $\breve {\boldsymbol e} =\mathscr L (\boldsymbol e)$ and $\breve
{\boldsymbol h} =\mathscr L (\boldsymbol h)$. Taking the Laplace transform of
\eqref{nrp} and eliminating $\breve {\boldsymbol h}$, we obtain
\begin{equation}\label{SDE}
\begin{cases}
\nabla \times \big( (\mu s)^{-1} \nabla \times \breve {\boldsymbol e}\big)
+\varepsilon s \breve {\boldsymbol e}=- \breve {\boldsymbol J} \quad
&{\rm in}~\Omega,\\
(\mu_j s)^{-1} (\nabla \times \breve {\boldsymbol e}) \times \boldsymbol n_j
+\mathscr B_j [\breve {\boldsymbol e}_{\Gamma_j}] =\breve {\boldsymbol V} \times
\boldsymbol n_j \quad &{\rm on}~\Gamma_j.
\end{cases}
\end{equation}

Our strategy is to show the well-posedness and stability of \eqref {SDE} in the
$s$-domain. The well-posedness of \eqref{nrp} follows from Lemma \ref{A2} and
the inverse Laplace transform. 

\begin{lemm}\label{es}
The problem (\ref{SDE}) has a unique weak solution $\breve {\boldsymbol e} \in
\boldsymbol H ({\rm curl}, \Omega)$ which satisfies
\begin{align}\label{st}
\|\nabla \times \breve {\boldsymbol e}\|_{\boldsymbol L^2 (\Omega)} &+\|s \breve
{\boldsymbol e}\|_{\boldsymbol L^2 (\Omega)}
\lesssim s_1 ^{-1}\Big[ \|s \breve {\boldsymbol J}\|_{\boldsymbol L^2
(\Omega)}\notag\\
&+ \sum\limits_{j=1}^2 ( \| s \breve {\boldsymbol V} \times \boldsymbol
n_j\|_{\boldsymbol H^{-1/2}({\rm div}, \Gamma_j)} +\| |s|^2 \breve {\boldsymbol
V} \times \boldsymbol n_j\|_{\boldsymbol H^{-1/2}({\rm div}, \Gamma_j)})\Big].
\end{align}
\end{lemm}

\begin{proof}
By Theorem \ref{AT}, it is easy to show the well-posedness of  the solution
$\breve {\boldsymbol e} \in \boldsymbol H ({\rm curl}, \Omega)$. 
Moreover, we have from the definition of \eqref{ASF} that 
\begin{align*}
a_{\rm TH} (\breve {\boldsymbol e}, \breve {\boldsymbol e})
&=\int_{\Omega} ( s \mu )^{-1} (\nabla \times \breve {\boldsymbol e}) \cdot
(\nabla \times \bar {\breve {\boldsymbol e}}){\rm d} \boldsymbol x
+\int_{\Omega} s \varepsilon \breve {\boldsymbol e} \cdot \bar {\breve
{\boldsymbol e}} {\rm d} \boldsymbol x +\sum\limits_{j=1}^2 \langle \mathscr B_j
[\breve {\boldsymbol e}_{\Gamma_j}], \breve {\boldsymbol e}_{\Gamma_j}
\rangle_{\Gamma_j}\\
&=-\int_{\Omega} \breve {\boldsymbol J} \cdot \bar {\breve {\boldsymbol e}} {\rm
d} \boldsymbol x+ \sum \limits_{j=1}^{2} \langle \breve {\boldsymbol V} \times
\boldsymbol n_j, \breve {\boldsymbol e}_{\Gamma_j} \rangle_{\Gamma_j}.
\end{align*}
It follows from the coercivity of $a_{\rm TH}$ in \eqref{RP} and the trace
theorem in Lemma \ref{TaT} that 
\begin{align*}
\frac{s_1}{|s|^2} \big(  \|\nabla \times \breve {\boldsymbol e}\|^2_{\boldsymbol
L^2 (\Omega)} +\|s \breve {\boldsymbol e}\|^2_{\boldsymbol L^2 (\Omega)} \big)
\lesssim &\|s^{-1} \breve {\boldsymbol J}\|_{\boldsymbol L^2 (\Omega)} \|s
\breve {\boldsymbol e}\|_{L^2 (\Omega)}+
\sum\limits_{j=1}^2 \|\breve {\boldsymbol V} \times \boldsymbol n_j\|_{\boldsymbol H^{-1/2} ({\rm div},  \Gamma_j)}
\| \breve {\boldsymbol e}_{\Gamma_j}\|_{\boldsymbol H^{-1/2} ({\rm curl}, \Gamma_j)}\\
\lesssim &\|s^{-1} \breve {\boldsymbol J}\|_{\boldsymbol L^2 (\Omega)} \|s
\breve {\boldsymbol e}\|_{L^2 (\Omega)}+
\sum\limits_{j=1}^2 \|\breve {\boldsymbol V} \times \boldsymbol n_j\|_{\boldsymbol H^{-1/2} ({\rm div},  \Gamma_j)}
\|\breve {\boldsymbol e}\|_{\boldsymbol H ({\rm curl}, \Omega)}\\
\lesssim &\|s^{-1} \breve {\boldsymbol J}\|_{\boldsymbol L^2 (\Omega)} \|s
\breve {\boldsymbol e}\|_{L^2 (\Omega)}
+\sum\limits_{j=1}^2 \|\breve {\boldsymbol V} \times \boldsymbol n_j\|_{\boldsymbol H^{-1/2} ({\rm div},  \Gamma_j)}
\|\nabla \times  \breve {\boldsymbol e} \|_{\boldsymbol L^2 (\Omega)}\\
&+\sum\limits_{j=1}^2 \|s ^{-1}\breve {\boldsymbol V} \times \boldsymbol n_j\|_{\boldsymbol H^{-1/2} ({\rm div},  \Gamma_j)}
\|s \breve {\boldsymbol e}\|_{\boldsymbol L^2 (\Omega)},
\end{align*}
which give the estimate \eqref{st} after applying the Cauchy--Schwarz
inequality.
\end{proof}

To show the well-posedness of the reduced problem \eqref {EHP}, we assume that 
\begin{equation}\label{H1}
 \boldsymbol E_0, \boldsymbol H_0 \in \boldsymbol H ({\rm curl }, \Omega), \quad \boldsymbol J \in
\boldsymbol H^1(0, T; \boldsymbol L^2 (\Omega)), \quad \boldsymbol J|_{t=0}=0.
\end{equation}
\begin{theo}
The problem \eqref{EHP} has a unique solution $(\boldsymbol E, \boldsymbol H)$, which satisfies
\begin{align*}
\boldsymbol E \in \boldsymbol L^2(0, T; \boldsymbol H ({\rm curl}, \Omega)) \cap
\boldsymbol H^1 (0, T; \boldsymbol L^2 (\Omega)),\\
\boldsymbol H \in \boldsymbol L^2(0, T; \boldsymbol H ({\rm curl}, \Omega)) \cap
\boldsymbol H^1 (0, T; \boldsymbol L^2 (\Omega)),
\end{align*}
and
\begin{align}
\int_0^{T} \big[\int_{\Omega} \big(\boldsymbol H \cdot  (\nabla \times
\bar{\boldsymbol \phi}) -\varepsilon \partial_t \boldsymbol E \cdot
\bar{\boldsymbol \phi}  \big){\rm d} \boldsymbol x- \sum \limits_{j=1}^2 \langle
\mathscr T_j \boldsymbol E_{\Gamma_j}, \boldsymbol \phi_{\Gamma_j}
\rangle_{\Gamma_j} \big] {\rm d} t \nonumber\\
=\int_0^T \int_{\Omega} \boldsymbol J \cdot \bar{\phi} {\rm d} \boldsymbol x
{\rm d} t, \quad \forall ~\boldsymbol \phi \in \boldsymbol H (\rm curl,
\Omega),\label{bf}\\
\int_0^T \int_\Omega (\nabla \times \boldsymbol E) \cdot \bar {\boldsymbol
\psi} +\mu \partial_t \boldsymbol H \cdot \bar {\boldsymbol \psi}
{\rm d} \boldsymbol x {\rm d} t=0, \quad \forall ~ \boldsymbol \psi \in
\boldsymbol L^2 (\Omega).
 \label{bf2}
\end{align}
Moreover, $(\boldsymbol E, \boldsymbol H)$ satisfy the stability estimate
\begin{align}
\max \limits_{[0, T]}
&\big( \|\partial_t \boldsymbol E\|_{\boldsymbol L^2 (\Omega)}+ \|\nabla \times
\boldsymbol E\|_{\boldsymbol L^2 (\Omega)} +\| \partial_t \boldsymbol
H\|_{\boldsymbol L^2 (\Omega)}+\|\nabla \times \boldsymbol H\|_{\boldsymbol L^2
(\Omega)}\big) \nonumber\\
&\lesssim\|\boldsymbol E_0\|_{\boldsymbol H (\rm curl, \Omega)}
+\|\boldsymbol{H}_0\|_{\boldsymbol H ({\rm curl}, \Omega)} +\|\boldsymbol
J\|_{\boldsymbol H^1
(0, T; \boldsymbol L^2 (\Omega))}.\label{ess}
\end{align}
\end{theo}

\begin{proof}
Let $\boldsymbol E= \boldsymbol U +\boldsymbol e$ and $\boldsymbol H=
\boldsymbol V +\boldsymbol h$, where $(\boldsymbol U,
\boldsymbol V)$ satisfy \eqref{TME} and $(\boldsymbol e, \boldsymbol h)$ satisfy
\eqref{nrp}. Noting 
\begin{align*}
&\int_0^T \big( \|\nabla \times \boldsymbol e\|_{\boldsymbol
L^2(\Omega)^2}^2+\|\partial_t \boldsymbol e\|^2_{\boldsymbol L^2(\Omega)}\big)
{\rm d} t\\
\leq & \int_0^T e^{- 2 s_1(t- T)}\big ( \|\nabla \times \boldsymbol e\|_{\boldsymbol L^2(\Omega)^2}^2
+\|\partial_t \boldsymbol e\|^2_{\boldsymbol L^2(\Omega)}\big) {\rm d} t\\
=& e^{2 s_1 T} \int_0^T e^{-2 s_1 t} \big ( \|\nabla \times \boldsymbol e\|_{\boldsymbol L^2(\Omega)^2}^2
+\|\partial_t \boldsymbol e\|^2_{\boldsymbol L^2(\Omega)}\big) {\rm d} t\\
\lesssim & \int_0^{\infty} e^{-2 s_1 t}  \big ( \|\nabla \times \boldsymbol e \|_{\boldsymbol L^2(\Omega)^2}^2
+\|\partial_t \boldsymbol e\|^2_{ \boldsymbol L^2(\Omega)}\big) {\rm d} t,
\end{align*}
we need to estimate
\[
\int_0^{\infty} e^{-2 s_1 t}  \big ( \|\nabla \times \boldsymbol e \|_{\boldsymbol L^2(\Omega)^2}^2
+\|\partial_t \boldsymbol e\|^2_{\boldsymbol L^2(\Omega)}\big) {\rm d} t.
\]
Taking the Laplace transform of (\ref {nrp}) yields
\begin{equation}\label{ln}
\begin{cases}
\nabla \times \breve {\boldsymbol e} +\mu s  \breve {\boldsymbol h}=0, \quad
\nabla \times \breve {\boldsymbol h} -\varepsilon s \breve {\boldsymbol
e}=\breve {\boldsymbol J}  \quad &{\rm in}~~\Omega,\\
\mathscr B_j [\breve {\boldsymbol e}_{\Gamma_j}]=\breve {\boldsymbol h} \times
\boldsymbol n_j+\breve {\boldsymbol V} \times n_j \quad
&{\rm on} ~\Gamma_j.
\end{cases}
\end{equation}
We have from Lemma \ref{es} that
\begin{align}\label{es1}
\|\nabla \times \breve {\boldsymbol e}\|_{\boldsymbol L^2 (\Omega)}+\|s
\breve {\boldsymbol e}\|_{\boldsymbol L^2 (\Omega)}\lesssim & s_1 ^{-1}
\Big[ \|s \breve {\boldsymbol J}\|_{\boldsymbol L^2
(\Omega)}+\sum\limits_{j=1}^2 ( \| s \breve {\boldsymbol V} \times \boldsymbol
n_j\|_{\boldsymbol H^{-1/2}({\rm div}, \Gamma_j)}\notag\\
 &+\| |s|^2 \breve {\boldsymbol V} \times \boldsymbol n_j\|_{\boldsymbol
H^{-1/2}({\rm div}, \Gamma_j)})\Big],
\end{align}
which gives after using (\ref {ln}) that 
\begin{align}\label{hs}
\|\nabla \times \breve {\boldsymbol h}\|_{\boldsymbol L^2 (\Omega)}+ \|s \breve
{\boldsymbol h}\|_{\boldsymbol L^2 (\Omega)}
&\lesssim s_1 ^{-1}
\big(\|\boldsymbol J\|_{\boldsymbol L^2 (\Omega)}+ \|s \breve {\boldsymbol
J}\|_{\boldsymbol L^2 (\Omega)}+ \sum\limits_{j=1}^2 ( \| s \breve {\boldsymbol
V} \times \boldsymbol n_j\|_{\boldsymbol H^{-1/2}({\rm div}, \Gamma_j)}\notag\\
 &+\| |s|^2 \breve {\boldsymbol V} \times \boldsymbol n_j\|_{\boldsymbol
H^{-1/2}({\rm div}, \Gamma_j)})\big).
\end{align}
It follows from \cite[Lemma 44.1]{Treves1975} that $\breve{\boldsymbol e}$ and
$\breve{\boldsymbol h}$ are holomorphic functions of $s$ on the half plane $s_1
>\bar\gamma>0,$  where $\bar \gamma$ is any positive constant. Hence we have
from Lemma \ref{A2} that the inverse Laplace transform of $\breve{\boldsymbol
e}$ and $\breve{\boldsymbol h}$ exist and are supported in $[0, \infty].$

Let $\boldsymbol e=\mathscr L^{-1} (\breve {\boldsymbol e})$ and $\boldsymbol
h =\mathscr L^{-1} (\breve {\boldsymbol h})$.
One may verify from the inverse Laplace transform and \eqref{a2} that $
\breve {\boldsymbol e} =\mathscr L (\boldsymbol e) =\mathscr F (e^{-s_1 t} 
\boldsymbol e)$, where $\mathscr F $ is the Fourier transform with respect to
$s_2$. It follows from the Parseval identity \eqref{PI} and \eqref{es1} that we
have
\begin{align*}
&\int_0^{\infty} e^{-2 s_1 t}  \big ( \|\nabla \times \boldsymbol e
\|_{\boldsymbol L^2(\Omega)^2}^2
+\|\partial_t \boldsymbol e\|^2_{\boldsymbol L^2(\Omega)}\big) {\rm d} t\\
&=\frac{1}{2 \pi} \int_{-\infty}^{\infty} \big(  \|\nabla \times
\breve{\boldsymbol e} \|_{\boldsymbol L^2(\Omega)^2}^2
+\|s\breve{\boldsymbol e}\|^2_{\boldsymbol L^2(\Omega)}\big) {\rm d} s_2\\
& \lesssim s_1 ^{-2} \int_{-\infty}^{\infty} \|s \breve {\boldsymbol
J}\|_{\boldsymbol L^2 (\Omega)} {\rm d} s_2+
s_1^{-2}\int_{-\infty}^{\infty} \sum\limits_{j=1}^2
 \big(\|s \breve {\boldsymbol V} \times \boldsymbol n_j\|^2_{\boldsymbol
H^{-1/2} ({\rm div}, \Gamma_j)}\\
&\hspace{2cm} + \||s|^2 \breve {\boldsymbol V} \times \boldsymbol
n_j\|^2_{\boldsymbol H^{-1/2} ({\rm div}, \Gamma_j)}
 \big){\rm d} s_2.
\end{align*}
By the assumption \eqref{H1}, we have $\boldsymbol J|_{t=0}=0$ in $\Omega$,
$\boldsymbol V \times \boldsymbol n_j|_{t=0}=\partial_t (\boldsymbol V
\times \boldsymbol n_j)|_{t=0}=0$ on $\Gamma_j$, which give that $\mathscr L
(\partial_t \boldsymbol J)= s \breve {\boldsymbol J}$ in $\Omega$ and $\mathscr
L (\partial_t (\boldsymbol V \times \boldsymbol n_j)) = s \breve {\boldsymbol V}
\times \boldsymbol n_j$ on $\Gamma_j.$ Noting
\[
|s|^2 \breve {\boldsymbol V} \times \boldsymbol n_j= (2 s_1-s) s\breve
{\boldsymbol V} \times \boldsymbol n_j = 2 s_1 \mathscr L (\partial_t
(\boldsymbol V\times \boldsymbol n_j)) -\mathscr L (\partial_t^2 (\boldsymbol V
\times \boldsymbol n_j)) \quad
{\rm on}~~ \Gamma_j,
\]
we have
\begin{align*}
&\int_0^{\infty}
e^{-2 s_1 t}  \big ( \|\nabla \times \boldsymbol e \|_{\boldsymbol
L^2(\Omega)^2}^2
+\|\partial_t \boldsymbol e\|^2_{\boldsymbol L^2(\Omega)}\big) {\rm d} t\\
&\lesssim s_1^{-2} \int_{-\infty}^{\infty} \|\mathscr L (\partial_t \boldsymbol
J)\|^2_{\boldsymbol L^2 (\Omega)} {\rm d} s_2 +s_1^{-2} \int_{-\infty}^{\infty} 
\sum\limits_{j=1}^2\|\mathscr L (\partial_t^2 (\boldsymbol V \times \boldsymbol
n_j))\|^2_{\boldsymbol H^{-1/2} ({\rm div}, \Gamma_j)}  {\rm d} s_2 \\
&\hspace{2cm}+(1+s_1^2)\int_{-\infty}^{\infty}\sum\limits_{j=1}^2 \| \mathscr L
(\partial_t (\boldsymbol V \times \boldsymbol n_j))\|^2_{\boldsymbol H^{-1/2}
({\rm div}, \Gamma_j)}  {\rm d} s_2.
\end{align*}
Using the Parseval identity \eqref{PI} again gives
\begin{align*}
&\int_0^{\infty}
e^{-2 s_1 t}  \big ( \|\nabla \times \boldsymbol e \|_{\boldsymbol
L^2(\Omega)^2}^2
+\|\partial_t \boldsymbol e\|^2_{\boldsymbol L^2(\Omega)}\big) {\rm d} t\\
&\lesssim s_1^{-2} \int_0^{\infty} e^{-2 s_1 t}\|\partial_t \boldsymbol
J\|^2_{\boldsymbol L^2 (\Omega)} {\rm d} t
+s_1^{-2} \int_0^{\infty} e^{-2 s_1 t} \sum\limits_{j=1}^2
\|\partial_t^2 (\boldsymbol V \times \boldsymbol n_j)\|^2_{\boldsymbol H^{-1/2}
({\rm div}, \Gamma_j)}  {\rm d} t\\
&\hspace{2cm}+(1+s_1^2)\int_0^{\infty}e^{-2 s_1 t} \sum\limits_{j=1}^2
\| \partial_t (\boldsymbol V \times \boldsymbol n_j)\|^2_{\boldsymbol H^{-1/2}
({\rm div}, \Gamma_j)}  {\rm d} t,
\end{align*}
which shows that
\[
\boldsymbol e \in \boldsymbol L^2 (0, T; \boldsymbol H ({\rm curl}, \Omega )
) \cap \boldsymbol H^{1} (0, T; \boldsymbol L^2 (\Omega)).
\]
Similarly, we can show from \eqref{hs} that
\[
\boldsymbol h \in \boldsymbol L^2 (0, T; \boldsymbol H ({\rm curl}, \Omega )
)\cap \boldsymbol H^{1} (0, T; \boldsymbol L^2 (\Omega)).
\]

Multiplying the test functions $\boldsymbol \psi \in \boldsymbol L^2 (\Omega) $
and $\boldsymbol \phi \in \boldsymbol H ({\rm curl}, \Omega)$  to the first and
second equality in \eqref{EHP}, respectively, using the boundary capacity
operators $\mathscr T_j$ and integration by parts,  we can get 
\eqref{bf}--\eqref{bf2}.

Next we show the stability estimate \eqref{ess}. Let $\tilde {
\boldsymbol E}$ be the extension of  $\boldsymbol E$ with respect to $t$ in
$\mathbb R$ such that $\tilde {\boldsymbol E} =0$ outside the interval $[0,
t]$. By the Parseval identity \eqref{PI} and  Lemma \ref{TP}, we get
\begin{align*}
 {\rm Re} \int_0^t \int_{\Gamma_j}e^{-2 s_1 t} \mathscr T_j [\boldsymbol
E_{\Gamma_j}] \cdot \bar {\boldsymbol E}_{\Gamma_j} {\rm d} \gamma_j {\rm d} t
 &= {\rm Re} \int_{\Gamma_j} \int_0^{\infty}
 e^{-2 s_1 t} \mathscr T_j [\tilde {\boldsymbol E}_{\Gamma_j}]\cdot \bar {
\tilde {\boldsymbol E}}_{\Gamma_j}{\rm d} \gamma_j
  {\rm d} t\\
 &=\frac{1}{2 \pi} \int_{-\infty}^{\infty} {\rm Re} \langle \mathscr B_j [\breve
{\tilde {\boldsymbol E}}_{\Gamma_j}],  \breve {\tilde{\boldsymbol
E}}_{\Gamma_j}\rangle_{\Gamma_j} {\rm d} s_2  \geq  0,
\end{align*}
which yields after taking $s_1 \rightarrow 0$ that
 \begin{equation}\label{TP1}
 {\rm Re} \int_0^t \int_{\Gamma_j} \mathscr T_j [\boldsymbol E_{\Gamma_j}] \cdot
\bar {\boldsymbol E}_{\Gamma_j} {\rm d} \gamma_j {\rm d} t \geq 0.
 \end{equation}

For any $0<t< T,$ consider the energy function
\[
e(t)= \|\varepsilon^{1/2} \boldsymbol E (\cdot, t)\|^2_{\boldsymbol L^2
(\Omega)} +\|\mu^{1/2} \boldsymbol H\|^2_{\boldsymbol L^2 (\Omega)}.
\]
It is  easy to note that
\[
\int_0^t e'(\tau) {\rm d} \tau =\big( \|\varepsilon^{1/2} \boldsymbol E (\cdot,
t) \|^2_{\boldsymbol L^2 (\Omega)} +\| \mu^{1/2} \boldsymbol H(\cdot,
t)\|^2_{\boldsymbol L^2 (\Omega)} \big)-\big( \|\varepsilon^{1/2} \boldsymbol
E_0\|^2_{\boldsymbol L^2 (\Omega)}
+\| \mu^{1/2} \boldsymbol H_0\|^2_{\boldsymbol L^2 (\Omega)} \big).
\]
On the other hand, it follows from \eqref{EHP}, \eqref{TP1}, and the integration
by parts that 
\begin{align}
\int_0^t e'(\tau) {\rm d} \tau
=& 2 {\rm Re} \int_0^t \int_{\Omega} (\varepsilon \partial_t \boldsymbol {E}
\cdot \bar {\boldsymbol E}+
\mu \partial_t \boldsymbol H \cdot \bar {\boldsymbol H}){\rm d} \boldsymbol x {\rm d} \tau \nonumber\\
=&2 {\rm Re}\int_0^t \int_{\Omega} \big((\nabla \times \boldsymbol H) \cdot
\bar {\boldsymbol E}-(\nabla \times \boldsymbol E) \cdot {\bar {\boldsymbol
H}}\big){\rm d} \boldsymbol x {\rm d} \tau -2 {\rm Re} \int_0^t \int_{\Omega}
\boldsymbol  J \cdot \bar {\boldsymbol E} {\rm d} \boldsymbol x {\rm d} \tau
\nonumber\\
=& 2 {\rm Re} \int_0^t \int_\Omega \big((\nabla \times \bar {\boldsymbol E})
\cdot \boldsymbol H- (\nabla \times \boldsymbol E) \cdot \bar {\boldsymbol H}
\big) {\rm d} \boldsymbol x {\rm d} \tau \nonumber\\
&-2{\rm Re} \sum \limits_{j=1}^2\int_0^t \int_{\Gamma_j} \mathscr
T_j [\boldsymbol E_{\Gamma_j}] \cdot \bar {\boldsymbol E}_{\Gamma_j} {\rm d}
\gamma_j {\rm d} \tau -2 {\rm Re} \int_0^t \int_{\Omega} \boldsymbol  J \cdot
\bar {\boldsymbol E} {\rm d} \boldsymbol x {\rm d} \tau \nonumber\\
\leq& -2 {\rm Re} \int_0^t \int_{\Omega} \boldsymbol  J \cdot \bar {\boldsymbol
E} {\rm d} \boldsymbol x {\rm d} \tau
\leq 2 \max \limits_{t \in [0,T]} \|\boldsymbol E\|_{\boldsymbol L^2
(\Omega)} \|\boldsymbol J\|_{\boldsymbol L^1(0, T; \boldsymbol L^2 (\Omega))}.
\label{es2}
\end{align}
Taking the derivative of \eqref{EHP} with respect to $t$, we know that
$(\partial_t \boldsymbol E, \partial_t \boldsymbol H)$ satisfy the
same set of equations  with the source $\boldsymbol J$ replaced by $\partial_t
\boldsymbol J$, and the initial conditions replaced by
$\partial_t \boldsymbol E|_{t=0}= \varepsilon^{-1 } \nabla \times
\boldsymbol H_0$, $\partial_t \boldsymbol H|_{t=0}= -\mu^{-1} \nabla \times
\boldsymbol E_0$. Hence we may follow the same steps as above to obtain
\eqref{es2} for $ (\partial_t \boldsymbol E, \partial_t \boldsymbol H),$ which
completes the proof of \eqref{ess} after combing the above estimates.
\end{proof}

\subsection{A priori estimates}

Now we intend to derive a priori stability estimates for the electric
field. Eliminating the magnetic field in \eqref{EH}--\eqref{IC} and using the
TBC in \eqref{ATTBC}, we consider the following initial-boundary value problem:
\begin{equation}\label{IB}
\begin{cases}
\varepsilon \partial_t^2 \boldsymbol E=-\nabla \times \big(  \mu^{-1} \nabla
\times \boldsymbol E\big ) -\boldsymbol F \quad &{\rm in}~\Omega, ~ t>0,\\
\boldsymbol E|_{t=0} =\boldsymbol E_0, \quad \partial_t \boldsymbol E|_{t=0}
=\boldsymbol E_1 \quad &{\rm in}~ \Omega,\\
\mu_j^{-1} (\nabla \times \boldsymbol E) \times \boldsymbol n_j +\mathscr  C_j
[\boldsymbol E_{\Gamma_j}]=0 \quad &{\rm on}~ \Gamma_j, ~ t>0,
\end{cases}
\end{equation}
where
\[
\boldsymbol F=\partial_t \boldsymbol J, \quad \boldsymbol E_1= \varepsilon^{-1}
(\nabla \times \boldsymbol H_0 -\boldsymbol J_0), \quad  \mathscr C_j = \mathscr
L^{-1} \circ s \mathscr B_j \circ \mathscr L.
\]
The variational problem \eqref{IB} is to find $\boldsymbol E \in \boldsymbol H
({\rm curl}, \Omega)$ for all $t>0$ such that
\begin{align}
\int_\Omega \varepsilon \partial_t^2 \boldsymbol E \cdot \bar {\boldsymbol w}
{\rm d } \boldsymbol x
=& -\int_{\Omega} \mu^{-1} (\nabla \times \boldsymbol E) \cdot (\nabla
\times \bar {\boldsymbol w} ) {\rm d} \boldsymbol x \nonumber\\
&-\int_{\Omega} \boldsymbol F \cdot \bar {\boldsymbol w} {\rm d} \boldsymbol x
-\sum\limits_{j=1}^2 \langle \mathscr C_j [\boldsymbol E_{\Gamma_j}],
\boldsymbol w_{\Gamma_j} \rangle_{\boldsymbol \Gamma_j},  \quad \forall ~
\boldsymbol w \in \boldsymbol H ({\rm curl }, \Omega).\label{tb}
\end{align}

\begin{lemm}\label{TTP}
Given $\xi \geq 0$ and $\boldsymbol E \in L^2 (0, \xi, \boldsymbol H^{-1/2}(\rm
curl, \Gamma_j))$, we have 
\[
{\rm Re} \int_0^{\xi} \int_{\Gamma_j}\left(  \int_0^t  \mathscr C_j [\boldsymbol
E_{\Gamma_j}](\tau) {\rm d} \tau \right) \cdot \bar{\boldsymbol E}_{\Gamma_j}
(t) {\rm d} \gamma_j {\rm d} t \geq 0.
\]
\end{lemm}
\begin{proof}
Let $\tilde {\boldsymbol E}$ be the extension of $\boldsymbol E$ with respect to
$t$ in $\mathbb R$ such that $\tilde {\boldsymbol E} =0$ outside the interval
$[0, \xi].$ It follows from the Parseval identity (\ref{PI}),  Lemma \ref{A2},
Lemma \ref{TP}, and \eqref{A1} that
\begin{align*}
  &{\rm Re }\int_{\Gamma_j}\int_{0}^{\xi} e^{-2 s_1 t} \big(  \int_0^t  \mathscr
 C_j [\boldsymbol E_{\Gamma_j}](\tau) {\rm d} \tau \big) \cdot \bar{\boldsymbol
E}_{\Gamma_j} (t) {\rm d} t {\rm d} \gamma_j \\
  &={\rm Re }\int_{\Gamma_j} \int_{0}^{\infty}e^{-2 s_1 t} \big(  \int_0^t 
\mathscr C_j [\tilde {\boldsymbol E}_{\Gamma_j}](\tau) {\rm d} \tau \big) \cdot
\bar { \tilde {\boldsymbol E}}_{\Gamma_j} (t)  { \rm d } t  {\rm d} \gamma_j\\
  &={\rm Re }\int_{\Gamma_j} \int_{0}^{\infty} e^{-2 s_1 t} \big( 
\int_0^t \mathscr  L^{-1} \circ  s \mathscr  B_j \circ \mathscr L \tilde
{\boldsymbol E}_{\Gamma_j} (\tau) {\rm d} \tau\big)\cdot \bar  {\tilde
{\boldsymbol E}}_{\Gamma_j} {\rm d } t  {\rm d} \gamma_j\\
  &=\frac{1}{2 \pi} {\rm Re } \int_{-\infty}^{\infty } \int_{\Gamma_j}  \mathscr
B_j \circ \mathscr L \tilde {\boldsymbol E}_{\Gamma_j}(s) \cdot \mathscr L (\bar
{\breve {\boldsymbol E}})(s ) {\rm d} \gamma_j {\rm d} s_2\\
  &=\frac{1}{2 \pi} \int_{-\infty}^{\infty} {\rm Re } \langle \mathscr
B_j[\breve {\tilde {\boldsymbol E}}_{\Gamma_j}], \breve {\tilde {\boldsymbol
E}}_{\Gamma_j} \rangle_{\Gamma_j} {\rm d }s_2 \geq 0.
\end{align*}
  The proof is completed by taking $s_1 \rightarrow  0$ in the above inequality.
\end{proof}
 
\begin{theo}\label{pe}
 Let $\boldsymbol E \in \boldsymbol H ({\rm curl}, \Omega)$ be the solution of
\eqref{tb}. If $\boldsymbol E_0, \boldsymbol E_1 \in \boldsymbol L^2 (\Omega)$
 and $\boldsymbol F \in \boldsymbol L^1 (0, T; \boldsymbol L^2 (\Omega)),$ 
 then $\boldsymbol E \in  \boldsymbol L^{\infty} (0, T; \boldsymbol L^2
(\Omega)).$
 Moreover, we have for any $T>0$ that 
 \begin{align}\label{IE}
\|\boldsymbol E\|_{\boldsymbol L^{\infty} (0, T; \boldsymbol L^2 (\Omega))}
\lesssim \|\boldsymbol E_0\|_{\boldsymbol L^2 (\Omega)} + T \| \boldsymbol
E_1\|_{\boldsymbol L^2 (\Omega)} + T \|\boldsymbol F\|_{\boldsymbol L^1 (0, T;
\boldsymbol L^2 (\Omega))},
 \end{align}
 and
 \begin{align}\label{IE2}
 \|\boldsymbol E\|_{\boldsymbol L^{2}(0, T; \boldsymbol L^2 (\Omega))} \lesssim
T^{1/2} \big(\|\boldsymbol E_0\|_{\boldsymbol L^2 (\Omega)} +
T \| \boldsymbol E_1\|_{\boldsymbol L^2 (\Omega)} + T \|\boldsymbol
F\|_{\boldsymbol L^1 (0, T; \boldsymbol L^2 (\Omega))} \big).
 \end{align}
 \end{theo}

 \begin{proof}
Let $ 0 < \xi < T $ and consider the function
\begin{equation}\label{af}
\boldsymbol \psi (\boldsymbol x, t) =\int_t^\xi
\boldsymbol E (\boldsymbol x, \tau) {\rm d} \tau, \quad \boldsymbol x \in
\Omega, ~ 0 \leq t \leq \xi.
 \end{equation}
 It is easy to verify that 
 \begin{equation}\label{p1}
 \boldsymbol \psi (\boldsymbol x, \xi)=0,\quad \partial_t \boldsymbol \psi
(\boldsymbol x, t) =- \boldsymbol E (\boldsymbol x, t),
 \end{equation}
 and
 \begin{equation}\label{p2}
 \int_0^\xi \boldsymbol \phi (\boldsymbol x, t) \bar {\boldsymbol \psi} (\boldsymbol x, t) {\rm d} t=\int_0^{\xi}
 \left( \int_0^t \boldsymbol \phi (\boldsymbol x, \tau) {\rm d} \tau  \right)
\cdot \bar {\boldsymbol E} (\boldsymbol x, t) {\rm d} t, \quad\forall ~
\boldsymbol \phi (\boldsymbol x, t) \in \boldsymbol L^2 (0, \xi; \boldsymbol L^2
(\Omega)). 
 \end{equation}
We show the last identity below. Using integration by parts and \eqref{p1} gives
 \begin{align*}
&\int_0^{\xi}\boldsymbol \phi(\boldsymbol x, t) \cdot  \bar { \boldsymbol \psi}
(\boldsymbol x, t) {\rm d} t =\int_0^{\xi } \left(\boldsymbol \phi(\boldsymbol
x, t) \cdot  \int_t^{\xi} \bar {\boldsymbol E}(\boldsymbol x, \tau) {\rm d} \tau
\right)  {\rm d} t\\
&=\int_0^{\xi} \int_t^{\xi} \bar  {\boldsymbol E} (\boldsymbol x,  \tau) {\rm d} \tau  \cdot  {\rm
d}\left(\int_0^t \boldsymbol  \phi(\boldsymbol x, \varsigma) {\rm d} \varsigma \right)\\
&=\int_t^{\xi} \bar {\boldsymbol E}(\boldsymbol x, \tau) {\rm d} \tau  \cdot
\int_0^t \boldsymbol \phi(\boldsymbol x, \varsigma) {\rm d} \varsigma
\Bigl{|}_0^{\xi}+\int_0^{\xi}\left( \int_0^t \boldsymbol \phi(\boldsymbol x,
\varsigma) {\rm d} \varsigma \right) \cdot \bar  {\boldsymbol E} (\boldsymbol
x, t) {\rm d} t\\
&=\int_0^{\xi} \left(\int_0^t \boldsymbol \phi(\boldsymbol x,  \tau) {\rm d}
\tau \right) \bar {\boldsymbol E}(\boldsymbol x, t)
{\rm d}  t.
\end{align*}

Taking the test function $\boldsymbol w= \boldsymbol \psi$ in (\ref{tb}) leads
to
\begin{align}
\int_\Omega \varepsilon \partial_t^2 \boldsymbol E \cdot \bar {\boldsymbol \psi}
{\rm d } \boldsymbol x
=& -\int_{\Omega} \mu^{-1} (\nabla \times \boldsymbol E) \cdot (\nabla
\times \bar {\boldsymbol \psi} ) {\rm d} \boldsymbol x \nonumber\\
&-\int_{\Omega} \boldsymbol F \cdot \bar {\boldsymbol \psi} {\rm d} \boldsymbol
x -\sum\limits_{j=1}^2 \langle \mathscr C_j [\boldsymbol
E_{\Gamma_j}], \boldsymbol \psi_{\Gamma_j} \rangle_{\boldsymbol \Gamma_j}.
\label{tf}
\end{align}
It follows from (\ref{p1})  and the initial conditions in (\ref{IB}) that
\begin{align*}
&{\rm Re} \int_0^\xi
 \int_\Omega \partial_t^2  \boldsymbol E \cdot \bar {\boldsymbol \psi} {\rm d}
\boldsymbol x {\rm d} t= {\rm Re} \int_{\Omega} \int_0^{\xi} \left( \partial_t
(\partial_t \boldsymbol E \cdot \bar {\boldsymbol \psi}) +\partial_t \boldsymbol
E \cdot
\bar {\boldsymbol E}\right) {\rm d} t{\rm d} \boldsymbol x\\
&={\rm Re}
\int_{\Omega} \left(  (\partial_t \boldsymbol E \cdot \bar {\boldsymbol
\psi})\Bigl |_{0}^{\xi} +\frac{1}{2} |\boldsymbol E|^2 \bigl |_0^\xi\right){\rm
d} \boldsymbol x\\
&=\frac{1}{2}\| \boldsymbol E (\cdot, \xi)\|^2_{\boldsymbol L^2
(\Omega)} -\frac{1}{2}\|\boldsymbol E_0\|^2_{\boldsymbol L^2 (\Omega)} - {\rm
Re} \int_{\Omega}\boldsymbol E_1 (\boldsymbol x) \cdot \bar {\boldsymbol \psi}
(\boldsymbol x, 0) {\rm d} \boldsymbol x.
\end{align*}
Thus, integrating (\ref{tf}) from $t=0$ to $t=\xi$ and taking the real parts
yields
\begin{align}
&\frac{\varepsilon}{2}\| \boldsymbol E (\cdot, \xi)\|^2_{\boldsymbol L^2
(\Omega)}
-\frac{\varepsilon}{2}\|\boldsymbol E_0\|^2_{\boldsymbol L^2 (\Omega)}
+\frac{1}{2}\int_{\Omega} \mu^{-1} \Bigl |  \int_0^ \xi \nabla \times
\boldsymbol E (\boldsymbol x, t) {\rm d} t \Bigl|^2 {\rm d} \boldsymbol x
\nonumber\\
&=\varepsilon {\rm Re} \int_{\Omega}
\boldsymbol E_1 (\boldsymbol x) \cdot \bar {\boldsymbol \psi} (\boldsymbol x, 0)
{\rm d} \boldsymbol x-{\rm Re} \int_0^\xi \int_{\Omega} \boldsymbol F \cdot \bar
{\boldsymbol \psi} {\rm d} \boldsymbol x {\rm d } t -{\rm Re}
\sum\limits_{j=1}^2\int_0^\xi  \langle \mathscr C_j [\boldsymbol E_{\Gamma_j}],
\boldsymbol \psi_{\Gamma_j} \rangle_{\boldsymbol \Gamma_j} {\rm d} t,\label{p4}
\end{align}
where we have used the fact that 
\[
\int_{\Omega} \mu^{-1} (\nabla \times \boldsymbol E) \cdot (\nabla \times \bar {\boldsymbol \psi} ) {\rm d} \boldsymbol x=
\frac{1}{2}\int_{\Omega} \mu^{-1} \Bigl |  \int_0^ \xi( \nabla \times
\boldsymbol E) {\rm d} t \Bigl|^2 {\rm d} \boldsymbol x.
\]
Next we estimate the three terms on the right-hand side of \eqref{p4}
separately.

We derive from  \eqref{af} and Cauchy--Schwarz inequality that
\begin{align}
&{\rm Re } \int_{\Omega} \boldsymbol E_1({\boldsymbol x}) \cdot \bar
{\boldsymbol \psi} (\boldsymbol x, 0) {\rm d} \boldsymbol x=
{\rm Re} \int_{\Omega} \boldsymbol E_1(\boldsymbol x) \cdot 
\left(  \int_0^{\xi} \bar {\boldsymbol E} (\boldsymbol x, t) {\rm d} t\right)
{\rm d} \boldsymbol x \nonumber\\
&={\rm Re}  \int_0^\xi \int_\Omega \boldsymbol E_1 (\boldsymbol x) \cdot \bar
{\boldsymbol E} (\boldsymbol x, t) {\rm d} \boldsymbol x {\rm d}t \leq
\|\boldsymbol E_1\|_{\boldsymbol L^2 (\Omega)} \int_0^\xi \|\boldsymbol E
(\cdot, t)\|_{\boldsymbol L^2 ( \Omega)}{\rm d} t. \label{I1}
\end{align}
Similarly, for $0 \leq t \leq \xi \leq T$, we have from \eqref{p2} that
\begin{align*}
&{\rm Re} \int_0^\xi \int_{\Omega} \boldsymbol F \cdot \bar {\boldsymbol \psi}
{\rm d} \boldsymbol x {\rm d} t
={\rm Re} \int_\Omega \int_0^\xi \left( \int_0^t \boldsymbol F (\boldsymbol x,
\tau) {\rm d} \tau\right) \cdot \bar {\boldsymbol E} (\boldsymbol x, t) {\rm d}
t {\rm d} \boldsymbol x \\
&={\rm Re} \int_0^\xi \int_0^t \int_{\Omega} \boldsymbol F (\boldsymbol x,
\tau)\cdot \bar {\boldsymbol E} (\boldsymbol x, t) {\rm d} \boldsymbol x{\rm d}
\tau {\rm d} t\\
&\leq \int_0^\xi \left(  \int_0^t \|\boldsymbol F (\cdot, \tau)\|_{\boldsymbol
L^2 (\Omega)} {\rm d} \tau\right) \|\boldsymbol E (\cdot, t)\|_{\boldsymbol
L^2(\Omega)} {\rm d}t  \\
&\leq \int_0^\xi \left(  \int_0^\xi \|\boldsymbol F (\cdot, \tau)\|_{\boldsymbol
L^2 (\Omega)} {\rm d} \tau\right) \|\boldsymbol E (\cdot, t)\|_{\boldsymbol
L^2(\Omega)} {\rm d}t\\
&\leq \left(  \int_0^\xi \|\boldsymbol F (\cdot, t)\|_{\boldsymbol L^2 (\Omega)}
{\rm d} t\right) \left (  \int_0^\xi \|\boldsymbol E (\cdot, t)\|_{\boldsymbol
L^2(\Omega)} {\rm d}t \right)
\end{align*}
Using Lemma \ref{TTP} and \eqref{p2}, we obtain
\begin{equation}\label{I3}
{\rm Re} \int_0^\xi  \langle \mathscr C_j [\boldsymbol E_{\Gamma_j}],
\boldsymbol \psi_{\Gamma_j} \rangle_{\boldsymbol \Gamma_j} {\rm d} t={\rm Re}
\int_0^{\xi} \int_{\Gamma_j}\left(  \int_0^t  \mathscr C_j [\boldsymbol
E_{\Gamma_j}](\tau) {\rm d} \tau \right) \cdot \bar{\boldsymbol E}_{\Gamma_j}
(t) {\rm d} \gamma_j {\rm d} t \geq 0.
\end{equation}
Substituting \eqref{I1}--\eqref{I3} into \eqref{p4}, we have for any $\xi \in
[0, T]$ that
\begin{align}
\frac{\varepsilon}{2}
&\| \boldsymbol E (\cdot, \xi)\|^2_{\boldsymbol L^2 (\Omega)}
+\frac{1}{2}\int_{\Omega} \mu^{-1}
\Bigl |  \int_0^ \xi \nabla \times \boldsymbol E (\boldsymbol x, t){\rm d} t
\Bigl|^2 {\rm d} \boldsymbol x \nonumber \\
&\leq\frac{\varepsilon}{2}\|\boldsymbol E_0\|^2_{\boldsymbol L^2 (\Omega)}
\left( \varepsilon \|\boldsymbol E_1\|_{\boldsymbol L^2 (\Omega)}+ \int_0^\xi
\|\boldsymbol F (\cdot, t)\|_{\boldsymbol L^2 (\Omega)} {\rm d} t\right) \left (
 \int_0^\xi \|\boldsymbol E (\cdot, t)\|_{\boldsymbol L^2(\Omega)} {\rm d}t
\right). \label{I4}
\end{align}
Taking the $\boldsymbol L^{\infty}$- norm with respect to $\xi$ on both sides of
\eqref{I4} yields
\[
\|\boldsymbol E\|^2_{\boldsymbol L^{\infty} (0, T; \boldsymbol L^2
(\Omega))}\lesssim \|\boldsymbol E_0\|^2_{\boldsymbol L^2 (\Omega)}+
T (\|\boldsymbol F\|_{\boldsymbol L^1 (0, T; \boldsymbol L^2
(\Omega))}+\|\boldsymbol E_1\|_{\boldsymbol L^2 (\Omega)}) \|\boldsymbol
E\|_{\boldsymbol L^{\infty} (0, T; \boldsymbol L^2 (\Omega))}.
\]
Therefore, the estimate \eqref{IE} follows directly from the Young
inequality. 

Integrating \eqref{I4} with respect to $\xi$ over $(0, T)$ and using
the Cauchy--Schwarz inequality, we obtain
\[
\| \boldsymbol E\|^2_{\boldsymbol L^2 (0, T; \boldsymbol L^2 (\Omega))}\lesssim
T\|\boldsymbol E_0\|^2_{\boldsymbol L^2 (\Omega)}+
T^{3/2} (\|\boldsymbol F\|_{\boldsymbol L^1 (0, T; \boldsymbol L^2
(\Omega))}+\|\boldsymbol E_1\|_{\boldsymbol L^2 (\Omega)}) \|\boldsymbol
E\|_{\boldsymbol L^{2} (0, T; \boldsymbol L^2 (\Omega))}.
\]
Using Young's inequality again, we derive the $\boldsymbol L^2$
estimate \eqref{IE2}, which completes the proof. 
\end{proof}

In Theorem \ref{pe}, it is required that $\boldsymbol E_0, \boldsymbol
E_1\in\boldsymbol L^2(\Omega)$, and $\boldsymbol F\in \boldsymbol L^1(0, T;
\boldsymbol L^2(\Omega))$, which can be satisfied if the data satisfy 
\[
 \boldsymbol E_0\in \boldsymbol L^2(\Omega),\quad\boldsymbol H_0\in\boldsymbol
H({\rm curl}, \Omega), \quad \boldsymbol J\in \boldsymbol H^1(0, T; \boldsymbol
L^2(\Omega)). 
\]

\section {Conclusion}\label{rem}

The scattering problems by unbounded structures have attracted much 
attention due to their wide applications and ample mathematical interests.
Although extensive study have been done for the time-harmonic problems, it is
still not clear what the best conditions are for those material parameters such
as the dielectric permittivity and magnetic permeability to assure the
well-posedness of the problems. In particular, it remains an open problem
whether it is well-posed for the real dielectric permittivity and magnetic
permeability. 

In this paper, we studied the time-domain scattering problem in an unbounded 
structure for the real dielectric permittivity and magnetic permeability. The
scattering problem was reduced to an initial-boundary value problem by using an
exact time-domain TBC. The reduced problem was shown to have a unique solution
by using the energy method. The main ingredients of the proofs were the Laplace
transform, the Lax--Milgram lemma, and the Parseval identity. Moreover, by
directly considering the variational problem of the time-domain wave equation,
we obtained a priori estimates with explicit dependence on time.

\end{document}